\documentclass[12pt]{amsart}

\usepackage{latexsym}
\usepackage{amsfonts}
\usepackage{amssymb}
\usepackage{amsmath}
\usepackage{mathrsfs}
\usepackage{hyperref}

\newcommand{\be}{\begin{equation}}
\newcommand{\ee}{\end{equation}}
\newcommand{\eea}{\end{eqnarray}}
\newcommand{\bean}{\begin{eqnarray*}}
\newcommand{\eean}{\end{eqnarray*}}
\newcommand{\brray}{\begin{array}}
\newcommand{\erray}{\end{array}}

\newtheorem{dfn}{Definition}[section]
\newtheorem{thm}[dfn]{Theorem}
\newtheorem{lmma}[dfn]{Lemma}
\newtheorem{ppsn}[dfn]{Proposition}
\newtheorem{cor}[dfn]{Corollary}
\newtheorem{xmpl}[dfn]{Example}
\newtheorem{rmrk}[dfn]{Remark}

\newcommand{\bdfn}{\begin{dfn}\rm}
\newcommand{\bthm}{\begin{thm}}
\newcommand{\blmma}{\begin{lmma}}
\newcommand{\bppsn}{\begin{ppsn}}
\newcommand{\bcor}{\begin{corcor}}
\newcommand{\bxmpl}{\begin{xmpl}}
\newcommand{\brmrk}{\begin{rmrk}\rm}

\newcommand{\edfn}{\end{dfn}}
\newcommand{\ethm}{\end{thm}}
\newcommand{\elmma}{\end{lmma}}
\newcommand{\eppsn}{\end{ppsn}}
\newcommand{\ecor}{\end{cor}}
\newcommand{\exmpl}{\end{xmpl}}
\newcommand{\ermrk}{\end{rmrk}}

\newcommand{\bbc}{\mathbb{C}}

\newcommand{\bbr}{\mathbb{R}}

\newcommand{\cla}{\mathcal{A}}
\newcommand{\clb}{\mathcal{B}}

\newcommand{\clh}{\ H }

\newcommand{\clk}{\mathcal{K}}

\newcommand{\e}{\mathrm{E}}

\newcommand{\en}{\mathrm{E}_0}

\newcommand{\m}{\ensuremath{\mathrm{M}}} 

\renewcommand{\k}{\ensuremath{\mathrm{k}}}

\title{$\e$-semigroups over closed convex cones}

\author{Anbu Arjunan}
\address{Chennai Mathematical Institute, H1, SIPCOT IT Park, Kelambakkam, Siruseri 603103, India.}
\email{aanbu@cmi.ac.in}

\author[R. Srinivasan]{R. Srinivasan} 
\email{vasanth@cmi.ac.in}
\thanks{Second author is currently visiting Kyoto University as a JSPS fellow}

\author[S. Sundar]{S. Sundar}
\email{ssundar@cmi.ac.in}


\begin{document}
\maketitle

\begin{abstract}
We initiate a study of $\e-$semigroups over convex cones.  We prove  a structure theorem for $\e-$semigroups which leave the algebra of compact operators invariant. Then we study in detail the CCR flows, $\en-$semigroups constructed from  isometric representations, by describing their units and gauge groups. We exhibit an uncountable family of $2-$parameter CCR flows, containing mutually non-cocycle-conjugate $\en-$semigroups.
\end{abstract}

\section{Introduction}
A one parameter $\en$-semigroup is a  $\sigma$-weakly continuous semigroup of normal unital $*-$endomorphisms on a von Neumann algebra, indexed by the positive real line. Since its inception in \cite{Pow}, the study of one-parameter $\en$-semigroups  on $B(H)$  has developed well in the last three decades. William Arveson contributed fundamentally to this development in a sequence of papers. In particular he classified completely the simplest class of $\en$-semigroups on $B(H)$ called as type I $\en$-semigroups. We refer to the monograph \cite{Arveson} for an extensive treatment regarding the theory of type I $\en-$semigroups on $B(H)$ (see also \cite{BVR}).

R.T  Powers discovered more complicated $\en$-semigroups  belonging to  type II and III (see \cite{powII} and \cite{powIII}.  In the last decade there were some significant  developments concerning these type II and III $\en$-semigroups on $B(H)$  (see \cite{T}, \cite{powII1}, \cite{pdct}, \cite{lieb}, \cite{gccr} and \cite{tcar}).  More recently there have been some important developments in the theory of $\en$-semigroups on non-type-I factors also. 

In this paper we initiate the theory of semigroups of endomorphisms indexed by a closed convex cone $P$ contained in $\bbr^d$. We  call them as $\en$-semigroups over $P$.  Already the relations between such $\en$-semigroups and the associated product systems of Hilbert spaces have been investigated (see \cite{Murugan_Sundar} and \cite{Murugan_Sundar 2}). Here, in this article, we discuss them systematically, with examples, and by computing their invariants like units, gauge groups, towards distinguishing them up to the equivalence of cocycle conjugacy. 

One parameter $\en$-semigroups arise naturally in the study of open quantum systems, the theory of interactions,  algebraic quantum field theory, and in quantum stochastic calculus.  Apart from their intrinsic interest,  arising as a natural mathematical  generalization,  we believe $\en$-semigroups over $P$ will have close connections with the study of quantum random fields and the theory of $C^*-$algebras. 

The simplest examples of  $\en-$semigroups are obtained by a process of second quantization from a semigroup of isometries. These are called as CCR flows. In the 1-parameter case, these CCR flows are of type I and they exhaust all type I examples. Type I is defined by the abundance of intertwining semigroup of isometries, that is such semigroups completely determine the $\en-$semigroup in some sense.  We can generalize the same process of second quantization to obtain CCR flows over any convex cone $P$. But in the multi-parameter case such abundance of units are not guaranteed  for CCR flows. Indeed all our examples admit only one unit, up to a multiple of scalars. This is the first complication we encounter in the multi-parameter case. 

On the other hand, Wold decomposition asserts that any strongly continuous one parameter semigroups of isometries is conjugate to a direct sum of semigroup of unitaries and a pure semigroup of isometries. Further any strongly continuous pure semigroup of isometries is conjugate to a right shift. There are countably many of them and they are determined by their index.  In the multi-parameter case, there is no Wold decomposition, and pure isometries are not determined uniquely by their index. The problem of classifying strongly continuous semigroup of isometries itself is still open. 

These certainly makes it hard to think about a possible classification of  `elementary $\en-$semigroups' like the type I examples in $1-$parameter case.  Indeed we exhibit  uncountably many non-conjugate semigroup of isometries, which lead to uncountably many non-cocycle-conjugate $\en-$semigroups over $\bbr_+\times \bbr_+$.  These are indications that there are more interesting mathematical objects and mathematical structures associated with multi-parameter $\en-$semigroups. We think  it is worth investigating them and this is a first step towards that.

 This paper is structured as follows. In Section 2, after fixing our notations we give the basic definitions of  $\en-$semigroup over convex cones, and the notion of  cocycle conjugacy. In Section 3, we recall the definition of multipliers and Wigner's theorem for automorphism groups. It follows that the automorphism groups are determined by the cohomology class of their multipliers.  We generalize a theorem of Arveson, which determine all $\e-$semigroups which leave the algebra of compact operators invariant, from the $1-$parameter case to the multi-parameter case. We call that as Arveson-Wigner's theorem.  In Section 4, we associate a CCR flow to an isometric representation of $P$, and provide the examples we are going to investigate in this paper. 
 
Unlike the $1-$parameter case, we need to take  cohomological considerations  in to account, while defining units for $\en-$semigroups over  convex cones.  In Section 5, after defining units,  we prove, for CCR flows, that the units are determined by the additive cocycles of  the associated isometric representations. We compute them for our specic examples.  In Section 6, we prove, for $\en-$semigroups in the so called standard form,  that cocycle conjugacy actually  implies conjugacy, a result known for $1-$parameter $\en$-semigroups.

In Section 7 and 8, we describe the gauge cocycles for CCR flows, and explicitly compute them for our specific examples. Using those computations, we exhibit a family containing uncountably many $2-$parameter $\en-$semigroups, which are mutually non-cocycle-conjugate.

\section{Preliminaries}
We fix the following notations which will be used throughout this paper.
Let $P \subset \mathbb{R}^{d}$ be a closed convex cone. We assume that $P$ is  spanning  and pointed, i.e. $P-P=\mathbb{R}^{d}$  and $P \cap -P=\{0\}$. Let $\Omega$ denote the interior of $P$. Then $\Omega$ is dense in $P$ (see Lemma 3.1 of \cite{Murugan_Sundar}).   
  Further  $\Omega$ is also spanning i.e. $\Omega-\Omega=\mathbb{R}^{d}$.  For $x, y \in \mathbb{R}^{d}$, we write $x \geq y$ and $x > y$ if $x-y \in P$ and $x-y \in \Omega$ respectively.

Our inner products are anti-linear in the first variable and linear in the second variable. Throughout this paper symbols like  $H, K$ will denote infinite dimensional separable complex Hilbert spaces, and the $*$-algebra of bounded operators on $H$ is denoted by  $B(H)$. We denote its predual, the ideal of trace class operators on $H$,  by $B^1(H)$, and the  algebra of compact operators on $H$  by $\mathcal{K}(H)$.  For $\xi,\eta \in H$, let $\theta_{\xi,\eta}$ be the rank-one operator on $H$ given by the equation $\theta_{\xi,\eta}(\gamma)=\xi\langle\eta|\gamma\rangle$.   For a von Neumann algebra $\m$, we denote by $\mathcal{U}(\m)$ the collection of all unitaries in $\m$, and for a Hilbert space $H$, we denote by $\mathcal{U}(H)$ the set of all unitary operators on $H$.
For Hilbert spaces $K, H$  and a unitary operator $U:H \to K$,  the map $B(H) \ni T \to UTU^{*} \in B(K)$ is denoted by $Ad(U)$. 

For a  complex separable Hilbert space $K$, we denote the symmetric Fock space by $\Gamma(K)$. We refer to \cite{KRP} for  proofs of the following well-known facts.  For $u \in K$, the exponential of $u$ is defined by  $e(u):=\sum_{n=0}^{\infty}\frac{u^{\otimes n}}{\sqrt{n!}}.$ Then the set $\{e(u): u \in K\}$ is linearly independent and total in $\Gamma(K)$. Exponential vectors satisfy $$\langle e(u),  e(v)\rangle = e^{\langle u | v\rangle} \;\;\text{for every}\; u,v\in K.$$   

For $u \in K$, there exists a  unitary operator, denoted $W(u)$, on $\Gamma(K)$ determined uniquely by the equation $W(u)e(v):=e^{-\frac{||u||^{2}}{2}-\langle u|v\rangle }e(u+v)\;\; \text{for}\; v\in K.$
 The operators $\{W(u): u \in K\}$ are called the Weyl operators and they  satisfy the following canonical commutation relation:  \[W(u)W(v)=e^{-i Im\langle u|v\rangle}W(u+v)~~ \forall u, v\in K.\]
 Further the linear span of $\{W(u): u \in K\}$ is a strongly dense unital $*$-subalgebra of $B(\Gamma(K))$. 
  For a unitary $U:K_1 \to K_2$ its second quantization is the unitary operator $\Gamma(U)$, from $\Gamma(K_1) \to \Gamma(K_2)$,  satisfying
$ \Gamma(U)e(v)=e(Uv)$ for $v \in K_1$.   
Second quantized unitaries are related to Weyl operators by  $$\Gamma(U)W(v)\Gamma(U)^{*}=W(Uv)~~ \forall v \in K_1.$$  
Second quantization can be defined for isometries as well in the same way. For two Hilbert spaces $K_1$ and $K_2$,  the map \[\Gamma(K_1 \oplus K_2) \ni e(u \oplus v) \to e(u) \otimes e(v) \in \Gamma(K_1) \otimes \Gamma(K_2)\] extends to a unitary operator. Through this unitary, we always identiy $\Gamma(K_1 \oplus K_2)$ with $\Gamma(K_1)\otimes \Gamma(K_2)$ without a mention.

For a complex separable Hilbert space $\k$ with finite or infinite dimension and $S\subseteq  \mathbb{R}^d$, $L^2(S, \k)$ denotes the Hilbert space of square integrable functions on the set $S$ taking values in $\k$, with respect to the Lebesgue measure. We end this section by recalling the basic definitions concerned with the  theory of $\en$-semigroups over $P$. 

\begin{dfn}
 An $\e$-semigroup over $P$ on a von Neumann algebra $\m$ is a family $\alpha = \{\alpha_x: x\in P\}$ of normal $*-$endomorphisms of $\m$ satisfying

\begin{itemize}
\item[(C1)] for $x, y \in P$, $\alpha_{x}\circ \alpha_{y}=\alpha_{x+y}$, $\alpha_0= Id_{\m}$ where $Id_{M}$ is the identity operator on $\m$ and
\item[(C2)] for $\rho\in \m_*$ and $m \in \m$, the map $P \ni x \to \rho(\alpha_x(m)) \in \mathbb{C}$ is continuous. 
\end{itemize}

 An $\e$-semigroup $\alpha$ is called  an $\en$-semigroup over $P$ if it is unital i.e, $\alpha_x(1)=1$ for all $x\in P$ and  said to be pure if $\cap_{t\geq 0} \alpha_{tx}(\m) =\bbc, ~~~\forall x\in \Omega.$  
\end{dfn}

In this paper we deal only with semigroups on $B(H)$. 
For $\alpha:=\{\alpha_{x}\}_{x \in P}$, a semigroup of normal $*$-endomorphisms on $B(H)$, 
using the fact that $*$-homomorphisms are contractive and the fact that finite rank operators are dense in $B^1(H)$, it is easy to see that Condition (C2) is equivalent to the following condition:
\begin{enumerate}
\item[(C2$^\prime$)] For $A \in B(H)$ and $\xi,\eta \in H$, the map $P \ni x \to \langle\alpha_{x}(A)\xi,\eta\rangle \in \mathbb{C}$ is continuous. 
\end{enumerate}
Since $P$ is fixed throughout, we simply refer an $\en$-semigroup over $P$ by an $\en$-semigroup.
We observe in the following lemma that the $\sigma-$weak continuity implies the strong continuity for $\mathrm{E}$-semigroups on $B(H)$.

\begin{lmma}
\label{strong continuity}
Let $\alpha:=\{\alpha_{x}\}_{x \in P}$ be an $\e$-semigroup on $B(H)$. Then  for $A \in B(H)$ and $\xi \in H$, the map $P \ni x \to \alpha_{x}(A)\xi \in H$ is continuous. 
\end{lmma}
\textit{Proof.}  
Let $(x_n)$ be a sequence in $P$ such that $(x_n) \to x$. Then
\begin{align*}
||\alpha_{x_n}(A)\xi-\alpha_{x}(A)\xi||^{2}
                                                               &=\langle \alpha_{x_n}(A^{*}A)\xi|\xi\rangle  - 2Re(\langle \alpha_{x_n}(A)\xi|\alpha_{x}(A)\xi\rangle)  +\langle \alpha_{x}(A^*A)\xi|\xi\rangle \\
                                                               &\to \langle \alpha_{x}(A^{*}A)\xi|\xi\rangle - 2 Re(\langle \alpha_{x}(A)\xi|\alpha_{x}(A)\xi\rangle +\langle \alpha_{x}(A^*A)\xi|\xi\rangle  \\
                                                               &=0.                                                
\end{align*}  \hfill $\Box$

\begin{dfn}
Let $\alpha:=\{\alpha_{x}\}_{x \in P}$ be an $\en$-semigroup on $B(H)$.
An $\alpha$-cocycle is  a strongly continuous family of unitaries $\{U_{x}\}_{x \in P}$ satisfying $$U_{x}\alpha_{x}(U_y)=U_{x+y}~~\forall x,y\in P.$$
A cocycle  $\{U_{x}\}_{x \in P}$ is said to be a gauge cocycle if further $U_x \in \alpha_x(B(H))^\prime$.
\end{dfn}

Given an $\alpha$-cocycle $\{U_{x}\}_{x \in P}$, it is easy to  verify that 
 $\{Ad(U_x)\circ \alpha_{x}\}$ is also an $\en$-semigroup on $B(H)$. 
 Let $\beta:=\{\beta_{x}\}_{x\in P}$ be another $\en$-semigroup on $B(H)$. We say $\beta$ is a cocycle perturbation of $\alpha$ if there exists an $\alpha$-cocycle
 $\{U_{x}\}_{x \in P}$ such that $\beta_{x}=Ad(U_x)\circ \alpha_{x}$ for every $x \in P$.
 
\begin{dfn}  Let $\alpha:=\{\alpha_{x}\}_{x \in P}$ and $\beta:=\{\beta_{x}\}_{x \in P}$ be $\en$-semigroups on $B(H)$ and $B(K)$ respectively. We say that

 \begin{itemize}
  \item[(i)]  $\alpha$ is  conjugate to $\beta$ if there exists a unitary operator $U:H \to K$ such that for every $x \in P$, $\beta_{x}=Ad(U)\circ \alpha_{x} \circ Ad(U^{*})$, and
   \item[(ii)]  $\alpha$ is cocycle conjugate to $\beta$ if there exists a unitary $U:H \to K$ such that $\en$-semigroup $\{Ad(U) \circ \alpha_{x} \circ Ad(U)^{*}\}_{x \in P}$ is a cocycle perturbation of $\beta$.  
\end{itemize}
\end{dfn}

Clearly cocycle conjugacy is an equivalence relation.

 \section{Arveson-Wigner's theorem}
 The main aim of this section is to establish a structure theorem for $\e$-semigroups which leave the algebra of compact operators invariant. We would like to call this theorem as Arveson-Wigner's theorem. We need some preparations before we can state the theorem.
 
  \begin{dfn}
 Let $\omega:P \times P \to \mathbb{T}$ be a continuous function. The map $\omega$ is called a multiplier on $P$ if 
 \[
 \omega(x,y)\omega(x+y,z)=\omega(x,y+z)\omega(y,z)\]
 for $x,y,z \in P$. 
  \end{dfn}
  The set of multipliers on $P$ forms an abelian group with respect to pointwise multiplication. 
  Let $\psi:P \to \mathbb{T}$ be a continuous function. Define \[\omega_{\psi}(x,y)=\psi(x)\psi(y)\psi(x+y)^{-1}.\] Then $\omega_{\psi}$ is a multiplier. Such multipliers are called coboundaries. Denote the set of multipliers on $P$ by $Z^{2}(P,\mathbb{T})$ and denote the set of coboundaries by $B^{2}(P,\mathbb{T})$. Note that $B^{2}(P,\mathbb{T})$ is a subgroup of $Z^{2}(P,\mathbb{T})$. Denote the quotient $Z^{2}(P,\mathbb{T})/B^{2}(P,\mathbb{T})$ by $H^{2}(P,\mathbb{T})$. For $\omega \in Z^{2}(P,\mathbb{T})$, $[\omega] $ denotes the class in $H^{2}(P,\mathbb{T})$ representing $\omega$.  Let $\omega_{1},\omega_2$ be multipliers on $P$. We say $\omega_1$ is cohomologous to $\omega_2$ if $[\omega_1]=[\omega_2]$.
  
\begin{xmpl} Let $A$ be an $d \times d$ real matrix. Define $\omega_{A}:P \times P \to \mathbb{T}$ by the formula
  \[
  \omega_{A}(x,y)=e^{i\langle Ax|y\rangle}
  \]
for $x,y\in P$.  Then $\omega_{A}(x,y)$ is a multiplier. 
\end{xmpl}

Any multiplier on $\bbr^d \times \bbr^d$ is equivalent to $\omega_A$ for some skew-symmetric $A$ (see Theorem 3.6.6 in  \cite{KRP1}). Also $\omega_{A_1}$ is equivalent to $\omega_{A_2}$ if and only if $A_1=A_2$. Now, using the original Wigner's theorem (for instance see Theorem 1.3.3 in  \cite{KRP1}), the isomorphism between skew-symmetric matrices and strictly upper triangular matrices, and the discussion in Section 1.4 in  \cite{KRP1}, the 
 following theorem can be deduced.

\begin{thm}[Wigner's theorem]
\label{Wigner}
Let $\{\alpha_{x}\}_{x \in \mathbb{R}^{d}}$ be a group of automorphisms on $B(H)$. Suppose that for $T \in B(H)$ and $\xi,\eta \in \clh$, the map $\mathbb{R}^{d} \ni x \to \langle\alpha_{x}(T)\xi|\eta\rangle \in \mathbb{C}$ is continuous. Then there exists a strongly continuous family of unitaries $\{U_{x}\}_{x \in \mathbb{R}^{d}}$ and a strictly upper triangular $d \times d$ real matrix $A$ such that 
\begin{enumerate}
\item[(1)] for $x \in \mathbb{R}^{d}$ and $T \in B(\clh)$, $\alpha_{x}(T)=U_{x}TU_{x}^{*}$, and
\item[(2)] for $x,y \in \mathbb{R}^{d}$, $U_{x}U_{y}=e^{i\langle A x|y\rangle}U_{x+y}$.
\end{enumerate}
\end{thm}

We need a few more lemmas.


\begin{lmma}
\label{Strong implies norm}
Let $p_{n}$ be a sequence of rank one projections in $B(H)$. Suppose that $p_{n}$ converges strongly to a rank one projection $p$. Then $p_n \to p$ in norm.
\end{lmma}
\textit{Proof.} Write $p_{n}=\theta_{\xi_n,\xi_n}$ and $p=\theta_{\xi,\xi}$ where $\xi_{n}$ and $\xi$ are unit vectors. Suppose $p_{n}$ does not converge to $p$ in norm. Then there exists $\epsilon > 0$ and a subsequence $p_{n_k}$ such that 
$||p_{n_k}-p|| \geq \epsilon$. By passing to a subsequence if necessary we can assume that $\xi_{n_k}$ converges weakly (say) to $\widetilde{\xi}$. 

We claim that there exists $\lambda \in \mathbb{C}$, with $|\lambda|=1$, such that $\widetilde{\xi}=\lambda \xi$. Let $\eta \in H$ be such that $\langle\xi|\eta\rangle=0$. Since $\theta_{\xi_{n_k},\xi_{n_k}}$ converges strongly to $\theta_{\xi,\xi}$, it follows that 
$\xi_{n_k}\langle\xi_{n_k}|\eta\rangle \to \xi \langle\xi|\eta\rangle=0$. Hence $|\langle\xi_{n_k}|\eta\rangle|=||\xi_{n_{k}}\langle\xi_{n_k}|\eta\rangle|| \to 0$. But $\xi_{n_k}$ converges weakly to $\widetilde{\xi}$. As a consequence, we have $\langle\widetilde{\xi}|\eta\rangle=0$.  This proves that 
there exists  $\lambda \in \mathbb{C}$ be such that $\widetilde{\xi}=\lambda \xi$.
Note that $\theta_{\xi_{n_k},\xi_{n_k}}(\xi)\to \theta_{\xi,\xi}(\xi)$.  This implies that \[|\langle\xi_{n_k}|\xi\rangle|=||\xi_{n_k}\langle\xi_{n_k}|\xi\rangle|| \to ||\theta_{\xi,\xi}(\xi)||=1.\]
But $\xi_{n_k}$ converges weakly to $\widetilde{\xi}=\lambda \xi$. Now the above convergence implies that $|\lambda|=1$. 

Note that $\theta_{\xi_{n_k},\xi_{n_k}}(\lambda \xi) \to \theta_{\xi,\xi}(\lambda\xi)$. Thus $\xi_{n_k}\langle\xi_{n_k}|\lambda \xi\rangle \to \xi\langle\xi|\lambda \xi\rangle=\lambda \xi$. But the sequence $\langle\xi_{n_k}|\lambda \xi\rangle \to \langle\lambda \xi|\lambda \xi\rangle=1$.  Calculate as follows to observe that 
\begin{align*}
\xi_{n_k}&=\xi_{n_k}\langle\xi_{n_k}|\lambda \xi\rangle\langle\xi_{n_k}|\lambda \xi\rangle^{-1} \\
              & \to \lambda \xi 
\end{align*}
Hence $\theta_{\xi_{n_k},\xi_{n_{k}}} \to \theta_{\lambda \xi,\lambda \xi}=\theta_{\xi,\xi}$ which is a contradiction. This implies that $p_n \to p$ in norm.  \hfill $\Box$

\begin{lmma}
\label{norm continuity}
Let $\alpha:=\{\alpha_{x}\}_{x \in P}$ be an $\e$-semigroup on $B(H)$. Suppose that for every $x \in P$, $\alpha_{x}$ leaves the algebra of compact operators invariant.  Then for every $T \in \mathcal{K}(H)$, the map $P \ni x \to \alpha_{x}(T) \in \mathcal{K}(H)$ is continuous where $\clk(H)$ is endowed with the norm topology.
\end{lmma}
\textit{Proof.} Fix $x \in P$. Then $\{\alpha_{tx}\}_{t \geq 0}$ is a $1$-parameter $\e$-semigroup on $B(H)$ leaving the algebra of compact operators invariant. By Prop. 3.4.1 of  \cite{Arveson}, it follows that there exists an isometry $V_{x}$ such that $\alpha_{x}(T)=V_{x}TV_{x}^{*}$. Consequently, it follows that for $x \in P$, $\alpha_{x}$ maps rank one projections to rank one projections. 

Since the linear combination of rank one projections is dense in $\clk(H)$, it suffices to show that if $p$ is a rank one projection then the map $P \ni x \to \alpha_{x}(p) \in \mathcal{K}(H)$ is norm continuous. Let $p$ be a rank one projection and let $(x_n)$ be a sequence in $P$ converging to $x \in P$. By Lemma \ref{strong continuity}, it follows that $\alpha_{x_n}(p) \to \alpha_{x}(p)$ strongly. Thanks to  Lemma \ref{Strong implies norm} we have $\alpha_{x_n}(p) \to \alpha_{x}(p)$   in norm. \hfill $\Box$

 \begin{rmrk}
 \label{Archimedean Property}
 Let $a \in \Omega$ and $x \in \mathbb{R}^{d}$ be given. Then there exists a positive integer $n$ such that $na>x$. This is due the fact that $a-\frac{x}{k} \to a \in \Omega$ as $k \to \infty$. Now the desired conclusion follows as $\Omega$ is an open convex cone. 
 \end{rmrk}

\begin{thm}[Arveson-Wigner's theorem]
\label{Arveson}
Let $\alpha:=\{\alpha_{x}\}_{x \in P}$ be an $\e$-semigroup on $B(H)$. Suppose that for every $x \in P$, $\alpha_{x}$ leaves the algebra of compact operators invariant. Then there exists a strongly continuous family $\{V_{x}\}_{x \in P}$ of isometries and a strictly upper triangular $d \times d$ real matrix such that 
\begin{enumerate}
\item[(1)] for $x \in P$, $\alpha_{x}(T)=V_{x}TV_{x}^{*}$, and
\item[(2)] for $x,y \in P$, $V_{x}V_{y}=e^{i\langle Ax,y\rangle}V_{x+y}$.
\end{enumerate}
\end{thm}
\textit{Proof.} We imitate the proof of Prop. 3.4.1 of \cite{Arveson}. We invoke Arveson's inductive limit construction. Fix a point $a \in \Omega$. Prop.3.4.1 of \cite{Arveson} applied to the $1$-parameter $\e$-semigroup $\{\alpha_{ta}\}_{t \geq 0}$ provides us with an isometry $V$ on $H$ such that $\alpha_{a}(T)=VTV^{*}$. Let $(U,\widetilde{H})$ be the minimal unitary dilation of $(V,H)$. This means the following.
\begin{enumerate}
\item[(1)] The Hilbert space $\widetilde{H}$ contains $H$ as a closed subspace. Strictly speaking, we have an isometric embedding of $H$ into $\widetilde{H}$. However to avoid cumbersome notations, we view $H$ as a closed subspace of $\widetilde{H}$. 
\item[(2)] The operator $U$ is a unitary operator on $\widetilde{H}$ such that $U\xi=V\xi$ for $\xi \in H$.
\item[(3)] The increasing union $\displaystyle \bigcup_{n\geq 0}U^{-n}H$ is dense in $\widetilde{H}$. 
\end{enumerate}
Decompose $\widetilde{H}$ as $H\oplus H^{\perp}$. The map $B(H) \ni T \to \begin{bmatrix}
                                                                                                                                        T & 0 \\
                                                                                                                                        0 & 0
                                                                                                                                        \end{bmatrix} \in B(\widetilde{H})$ is a $*$-homomorphism and is an embedding. Thus we view $B(H)$ as a $*$-subalgebra of $B(\widetilde{H})$. 
                                                                                                                                        
 For $n\geq 0$, let $\clk_{n}$ be the closed linear span of $\{\theta_{\xi,\eta}: \xi,\eta \in U^{-n}H\}$.    Note that $\clk_{n}$ is an increasing union of $C*$-subalgebras of $\mathcal{K}(\widetilde{H})$. Also the union $\displaystyle \bigcup_{n \geq 0} \clk_{n}$ is dense in $\mathcal{K}(\widetilde{H})$.   Note that $\clk_{0}=\clk(H)$.    We have the following.
 \begin{enumerate}
 \item[(1)] For $T \in \clk_{0}$, $\alpha_{a}(T)=UTU^{*}$. 
 \item[(2)] Let $n \geq 0$ be given. Note that the map $\clk_{n} \ni T \to U^{n}TU^{-n} \in \clk_{0}$ is an isomorphism. For $x \in P$, let $\beta_{x}^{(n)}:\clk_{n} \to \clk_{n}$ be defined by the equation \[\beta_{x}^{(n)}(T)=U^{-n}\alpha_{x}(U^{n}TU^{-n})U^{n}.\]
 It is clear that $\beta_{x}^{(n)}$ is a $*$-homomorphism.
 Note that for $x,y \in P$, $\beta_{x}^{(n)}\beta_{y}^{(n)}=\beta_{x+y}^{(n)}$. For $T \in K_{n}$, $\beta_{a}^{(n)}(T)=UTU^{*}$. Fix $T \in \clk_{n}$. Lemma \ref{norm continuity} implies that the map $P \ni x \to \beta_{x}^{(n)}(T) \in \clk_{n}$ is continuous.
 \item[(3)] Let $n \geq 0$ and $T \in \clk_{n}$ be given. Using $(1)$, it is routine to check that for $x \in P$, \[\beta_{x}^{(n+1)}(T)=\beta_{x}^{(n)}(T).\] 
 Using the fact that $\clk_{n}$ is simple for every $n \geq 0$, we deduce that for $x \in P$, there exists a $*$-homomorphism, denoted $\beta_{x}$, of $\clk(\widetilde{H})$ such that $\beta_{x}(T)=\beta_{x}^{(n)}(T)$ for every $T \in \clk_{n}$. Observe that $\beta_{a}(T)=UTU^{*}$.  Also note that for $x \in P$ and $T \in \clk_{0}$, $\beta_{x}(T)=\alpha_{x}(T)$. 
 
 \item[(4)] It is clear from $(2)$ that for $T \in \clk(\widetilde{H})$, the map $P \ni x \to \beta_{x}(T) \in \mathcal{K}(\widetilde{H})$ is continuous. Here $\clk(\widetilde{H})$ is given the norm topology.
 
 \item[(5)] From $(2)$, it is clear that $\{\beta_{x}\}_{x \in P}$ is a semigroup of $*$-endomorphisms of $\clk(\widetilde{H})$. Let $x \in P$ be given. We claim that  $\beta_{x}$ is an automorphism. Since $\clk(\widetilde{H})$ is simple, it suffices to show that $\beta_{x}$ is onto. For $y \in P$, let $A_{y}:=\beta_{y}(\clk(\widetilde{H}))$. Note that $(A_{y})_{y \in P}$ is decreasing with respect to the order $<$. Since $\beta_{a}$ is an automorphism, it follows that $A_{na}=\clk(\widetilde{H})$ for every $n$. Now the desired conclusion follows from Remark \ref{Archimedean Property}.
 \item[(6)] For $x \in P$, we denote the extension of $\beta_{x}$ to the multiplier algebra of $\clk(\widetilde{H})$, which is $B(\widetilde{H})$, by $\beta_{x}$ itself. We claim that for $T \in B(\widetilde{H})$ and $\xi,\eta \in \widetilde{H}$, the map $P \ni x \to \langle\beta_{x}(T)\xi|\eta\rangle \in \mathbb{C}$ is continuous.  Let $T \in B(\widetilde{H})$ and $\xi,\eta \in \widetilde{H}$ be given.  If $T \in \clk(\widetilde{H})$, then by $(4)$, it follows that the map $P \ni x \to \langle\beta_{x}(T)\xi|\eta\rangle \in \mathbb{C}$ is continuous. Otherwise, let $(T_n)$ be a sequence in $\clk(\widetilde{H})$ such that $T_n \to T$ strongly and $||T_{n}|| \leq ||T||$. Then $\beta_{x}(T_n) \to \beta_{x}(T)$ weakly for every $x \in P$.  This implies that the sequence of continuous functions $\langle\beta_{x}(T_n)\xi|\eta\rangle$ converges pointwise to $\langle\beta_{x}(T)\xi|\eta\rangle$. Hence the map $P \ni x \to \langle\beta_{x}(T)\xi|\eta\rangle$ is measurable for every $T \in B(\widetilde{H})$ and $\xi,\eta \in \widetilde{H}$. By Corollary 4.3 of \cite{Murugan_Sundar}, we have for $T \in B(\widetilde{H})$ and $\xi,\eta \in \widetilde{H}$, the map $P \ni x \to \langle\beta_{x}(T)\xi|\eta\rangle$ is continuous. This proves our claim. 
 
 \item[(7)] For $z \in \mathbb{R}^{d}$, write $z=x-y$ with $x,y \in P$. Set $\widetilde{\beta}_{z}=\beta_{x}\circ \beta_{y}^{-1}$. We leave it to the reader to verify that $\widetilde{\beta}_{z}$ is well-defined. It is clear that $\{\widetilde{\beta}_{z}\}_{z \in \mathbb{R}^{d}}$ is a group of automorphisms of $B(\widetilde{H})$. Note that for $x \in P$, $\widetilde{\beta}_{x}=\beta_{x}$. Let $T \in B(\widetilde{H})$ and $\xi,\eta \in \widetilde{H}$ be given. Let $z_n$ be a sequence in $\mathbb{R}^{d}$ such that $z_n \to z \in \mathbb{R}^{d}$. Since $\mathbb{R}^{d}=\Omega-\Omega$, write $z=x-y$ with $x,y \in \Omega$. Since $z \in \Omega-y$, it follows that $z_n \in \Omega-y$ eventually. Thus for $n$ large, there exists $x_n \in \Omega$ such that $z_n=x_n-y$. Now $x_n \to x$. Now calculate as follows to observe that 
 \begin{align*}
 \langle\widetilde{\beta_{z_n}}(T)\xi|\eta\rangle&=\langle\beta_{x_n}(\beta_{y}^{-1}(T))\xi|\eta\rangle \\
                                                        & \to \langle\beta_{x}(\beta_{y}^{-1}(T))\xi|\eta\rangle ~~(\textrm{by (6)})\\
                                                         & = \langle\widetilde{\beta}_{z}(T)\xi|\eta\rangle
  \end{align*}
 This proves that for $T \in B(\widetilde{H})$ and $\xi,\eta \in \widetilde{H}$, the map $\mathbb{R}^{d} \ni z \to \langle\widetilde{\beta}_{z}(T)\xi|\eta\rangle \in \mathbb{C}$ is continuous.  By Wigner's theorem (\ref{Wigner}), there exists a strongly continuous family of unitaries $\{W_{z}\}_{z \in \mathbb{R}^{d}}$ and a strictly upper triangular matrix $d \times d$ matrix $A$ such that 
 \begin{align*}
 \widetilde{\beta}_{z}(T)&=W_{z}TW_{z}^{*} \\
 W_{z_1}W_{z_2}&=e^{i\langle Az_1|z_2\rangle}W_{z_1+z_2}
  \end{align*}
 for $T \in B(\widetilde{H})$ and $z,z_1,z_2 \in \mathbb{R}^{d}$.
 
 \item[(8)] For $x \in P$ and $T \in \clk_{0}$, $W_{x}TW_{x}^{*}=\alpha_{x}(T) \in \clk_{0}$. Consequently, for $x \in P$, $W_{x}$ leaves $H$ invariant. For $x \in P$, let $V_{x}:H \to H$ be defined by $V_{x}\xi=W_{x}\xi$ for $\xi \in H$. 
 It is clear that for $x,y \in P$, \[V_{x}V_{y}=e^{i\langle Ax|y\rangle}V_{x+y}.\] We leave it to the reader to verify that 
  for $x \in P$ and $T \in \mathcal{K}(H)$, $\alpha_{x}(T)=V_{x}TV_{x}^{*}$.   Since $\clk(H)$ is dense in $B(H)$ in the $\sigma$-weak topology, it follows that for $x \in P$ and $T \in B(H)$, $\alpha_{x}(T)=V_{x}TV_{x}^{*}$.
  \end{enumerate}
 This completes the proof. \hfill $\Box$
 
 Let $B_{0}(d)$ be the additive group of strictly upper triangular real matrices. For $A \in B_{0}(d)$, let $\omega_{A}$ be the multiplier on $P$ defined by $\omega_{A}(x,y)=e^{i\langle Ax|y\rangle}$.
The following result is known in the measurable setting \cite{Laca}
\begin{cor}
\label{second cohomology}
The map $B_{0}(d) \ni A \to [\omega_{A}] \in H^{2}(P,\mathbb{T})$ is an isomorphism.
\end{cor}
\textit{Proof.} It is evident that the map $B_{0}(d) \ni A \to [\omega_{A}]\in H^{2}(P,\mathbb{T})$ is a homomorphism. Let $A \in B_{0}(d)$ be given. Suppose $\omega_{A}$ is a coboundary. Then $\omega_{A}(tx,y)=\omega_{A}(y,tx)$ for $x,y \in A$ and $t >0$. This implies that for $t >  0$ and $x,y \in P$, $e^{it\langle Ax|y\rangle}=e^{it\langle Ay|x\rangle}$. Hence $\langle Ax|y\rangle=\langle Ay|x\rangle$ for all $x,y\in P$. Since $P$ is spanning, it follows that $\langle Ax|y\rangle=\langle Ay|x\rangle$ for all $x,y \in \mathbb{R}^{d}$. This implies $A=A^{t}$. But $A$ is strictly upper triangular. Hence $A=0$. This shows that the map $B_{0}(d) \ni A \to [\omega_{A}] \in H^{2}(P,\mathbb{T})$ is one-one. 

Let $\omega$ be a multiplier on $P$. Consider the Hilbert space $H:=L^{2}(P)$. For $x \in P$, let $V_{x}$ be the isometry on $H$ defined by the formula: \begin{equation}\label{standard}
V_{x}(f)(y):=\begin{cases}
 \omega(x,y-x)f(y-x)  & \mbox{ if
} y -x \in P,\cr
   &\cr
    0 &  \mbox{ if } y-x \notin P.
         \end{cases}
\end{equation}
for $f \in H$.  Then $\{V_{x}\}_{x \in P}$ is a strictly continuous family of isometries  i.e. for $f \in H$, the maps $P \ni x \to V_{x}f\in H$ and $P\ni x \to V_{x}^{*}f \in H$ are continuous. It is clear that for $x,y \in P$, $V_{x}V_{y}=\omega(x,y)V_{x+y}$.

 For $x \in P$, let $\alpha_{x}$ be the endomorphism of $B(H)$ given by the formula 
 $
 \alpha_{x}(T)=V_{x}TV_{x}^{*}
 $
 for $T \in B(H)$. Then $\alpha:=\{\alpha_{x}\}_{x \in P}$ is an $\e$-semigroup leaving the algebra of compact operators invariant. Thus by Arveson-Wigner's theorem (\ref{Arveson}), it follows that there exists a strongly continuous family $\{W_{x}\}_{x \in P}$ of isometries and a strictly upper triangular $d\times d$ real matrix $A$ such that 
 \begin{enumerate}
 \item[(1)] for $T \in B(H)$ and $x \in P$, $\alpha_{x}(T)=W_{x}TW_{x}^{*}$, and
\item[(2)] for $x,y \in P$, $W_{x}W_{y}=\omega_{A}(x,y)W_{x+y}$.
\end{enumerate}
Fix $x \in P$. Let $T \in B(H)$ be given. Calculate as follows to observe that 
\begin{align*}
V_{x}^{*}W_{x}T&=V_{x}^{*}W_{x}TW_{x}^{*}W_{x} \\
                          &=V_{x}^{*}V_{x}TV_{x}^{*}W_{x}\\
                          &=TV_{x}^{*}W_{x}.
\end{align*}
Thus $V_{x}^{*}W_{x}$ commutes with every element of $B(H)$. It follows that $V_{x}^{*}W_{x}$ is a scalar which we denote by $f(x)$. The strong continuity of $\{V_{x}\}_{x\in P}$ and $\{W_{x}\}_{x \in P}$ implies that $f$ is continuous. Now calculate as follows to observe that for $x \in P$, 
\[W_{x}=W_{x}W_{x}^{*}W_{x}=V_{x}V_{x}^{*}W_{x}=f(x)V_{x}.\] This implies that $f$ takes values in $\mathbb{T}$. Calculate as follows to observe that for $x,y \in P$, 
\begin{align*}
\omega_{A}(x,y)&=W_{x+y}^{*}W_{x}W_{y} \\
                          &=f(x+y)^{-1}f(x)f(y)V_{x+y}^{*}V_{x}V_{y}\\
                          &=f(x+y)^{-1}f(x)f(y)\omega(x,y).
\end{align*}
This shows that $[\omega]=[\omega_{A}]$. This completes the proof. \hfill $\Box$

\section{Examples}

\begin{dfn} Let $K$ be a separable Hilbert space.  A map $V:P \to B(K)$ is called  an isometric representation of $P$ on $K$ if
\begin{enumerate}
\item[(1)]  $V_{x}$ is an isometry for all $x \in P$,
\item[(2)]  $V_{x+y}=V_{x}V_{y}$ for all $x,y \in P$, and
\item[(3)] the map $P \ni x \to V_{x}\xi \in K$ is continuous  for all  $\xi \in K$. 
\end{enumerate}

 An isometric representation is said to be pure if $\cap_{t\geq 0} V_{ta} (K) = \{0\}$ for all $a \in \Omega$.  $V$ is strictly continuous if  further
the map $P \ni x \to V_{x}^{*}\xi \in K$ is also continuous  for all $\xi \in K$.  
\end{dfn}
There is an obvious addition operation on the class of isometric representations of $P$. Let $V:=\{V_{x}\}_{x \in P}$ and $W:=\{W_{x}\}_{x \in P}$ be isometric representations of $P$ on the Hilbert spaces $H$ and $K$ respectively. Then $V\oplus W:=\Big\{\begin{bmatrix}
 V_{x} & 0 \\
 0  & W_{x}
 \end{bmatrix}\Big\}_{x \in P}$ is an isometric representation of $P$ on $H \oplus K$. Clearly $V \oplus W$ is strictly continuous if $V$ and $W$ are strictly continuous.  
 Here is a basic example of an isometric representation.
 
 \begin{xmpl}\label{usualshift}
 Let $\k$ be a separable Hilbert space.  For $x \in P$, let $S_{x}$ be the shift operator on $K = L^2(P, \k)$ defined by \begin{equation}
S_{x}(f)(y):=\begin{cases}
 f(y-x)  & \mbox{ if
} y -x \in P,\cr
   &\cr
    0 &  \mbox{ if } y-x \notin P.
         \end{cases}
\end{equation}
 Then $S =\{S_{x}\}_{x \in P}$ is a strictly continuous isometric representation. Throughout this paper,  we refer to this representation by $\{S_{x}\}_{x \in P}$. 
 \end{xmpl}
 
 Unlike in the $1$-parameter case (when $P = \bbr_+$), where  all pure isometric representations are  conjugate to the right shift, there are many non-conjugate  isometric representations in the multi-parameter case.  
 
 \begin{xmpl}
Let $A \subseteq \mathbb{R}^{d}$ be a nonempty closed subset. We say that $A$ is a $P$-module if $A+P \subseteq A$.  
Consider the Hilbert space $K:= L^{2}(A,\k)$. For $x \in P$, let $S^A_{x}$ be the operator on $K$ defined by \begin{equation}
S^A_{x}(f)(y):=\begin{cases}
 f(y-x)  & \mbox{ if
} y -x \in A,\cr
   &\cr
    0 &  \mbox{ if } y-x \notin A.
         \end{cases}
\end{equation}
 Then $S^A= \{S^A_{x}\}_{x \in P}$ is a strictly continuous isometric representation. We call $\{S^A_{x}\}_{x \in P}$ the isometric representation associated to the $P$-module $A$ of multiplicity $\dim(\k)$.
  \end{xmpl}

 \begin{xmpl}\label{A12}
 When $P = \bbr_+\times \bbr_+$,
it is easy to check $(S_{(0,t)})^*S_{(s,0)} = S_{(s,0)} (S_{(0,t)})^*$ for all $(s,t) \in \bbr_+\times \bbr_+$. But this relation is not satisfied by $S^{A_T^{(a,b)}}$, where $$A_T^{(a,b)} =  [a, \infty) \times [b, \infty) ~\bigsqcup ~ [0,\infty) \times [0,b), ~~~ a<0, ~0<b.$$ Indeed, for $ (u,v)\in (0,\infty) \times (0,b)$ and $ u<s < u+a,  b-v< t$, $((S^{A_T^{(a,b)}}_{(0,t)})^*S^{A_T^{(a,b)}}_{(s,0)}f) (u,v) =0$, but $(S^{A_T^{(a,b)}}_{(s,0)} (S^{A_T^{(a,b)}}_{(0,t)})^*f)(u,v)$ will not be equal to $0$, for a suitable choice of $f$, which we leave as an easy exercise to the reader. This shows that the isometric representations $S$ and $S^{A_T^{(a,b)}}$ are not conjugate.
 \end{xmpl}

Let $A$ be a $P$-module. Then $A+P \subset A$. Hence by Lemma II.12 of \cite{Hilgert_Neeb}, $Int(A)$ is dense in $A$ and the boundary of $A$, $\partial A$ has 
zero measure.

\begin{lmma}
\label{translate1}
Let $A \subset \mathbb{R}^{d}$ be a $P$-module and  $a \in \Omega$ be given. Assume that $A \neq \mathbb{R}^{d}$. Then there exists $n \geq 1$ such that $(A+na) \cap -\Omega =\emptyset$. 
\end{lmma}
\textit{Proof.} Let $A \subset \mathbb{R}^{d}$ and $a \in \Omega$ be as in the statement of the lemma. Suppose for every $n \geq 1$, $(A+na)\cap -\Omega \neq \emptyset$. Then for $n \geq 1$, there exists $x_n \in A$ and $y_n \in \Omega$ such that 
$x_n+na=-y_n$. This implies that $-na=x_n+y_n \in A$ for every $n$. For $A + \Omega \subset A$. Let $x \in \mathbb{R}^{d}$ be given. By Remark  \ref{Archimedean Property}, there exists a positive integer $n_0$ such that  $n_0a+x \in \Omega$.  This implies that $x=-n_0a+(x+n_0a) \in A$.  Consequently, we have $A=\mathbb{R}^{d}$ which is a contradiction. This contradiction proves that there exists $n \geq 1$ such that $(A+na)\cap -\Omega=\emptyset$. This completes the proof. \hfill $\Box$

\begin{ppsn}
\label{purity}
Let $A \subset \mathbb{R}^{d}$ be a $P$-module. Assume that $A\neq \mathbb{R}^{d}$.  Then the isometric representation associated to $A$ of any multiplicity is pure.
\end{ppsn}
\textit{Proof.} Let $V:=\{V_{x}\}$ be the isometric representation associated to $A$ of multiplicity $k$. We can assume that $k=1$. Since $V$ is unitarily equivalent to the isometric representation associated to a module which is  a translate of $A$, by Lemma \ref{translate1} we can assume that $A \cap -\Omega=\emptyset$. We make both the assumptions. Since $A$ and $Int(A)$ differ by a set of a measure zero, it follows that $C_{c}(Int(A))$ is dense in $L^{2}(A)=L^{2}(Int(A))$. 

Let $\phi \in C_{c}(Int(A))$ and $a \in \Omega$ be given. It is enough to show that $V_{ta}^{*}\phi \to 0$ as $t \to \infty$.  Denote the support of $\phi$ by $K$. We claim that there exists $t_{0} \geq 0$ such that $t \geq t_0$ implies that $K \cap (A+ta)=\emptyset$. 
Suppose not. Then there exists a sequence $(t_n) \to \infty$ and $x_{n} \in K$ and $y_n \in A$ such that $x_n=y_n+t_na$. Note that $(x_n)$ is bounded. Hence $\frac{x_n}{t_n} \to 0$. This implies that $\frac{y_n}{t_n} \to -a \in -\Omega$. Hence eventually $\frac{y_n}{t_n} \in -\Omega$. Thus eventually $y_{n} \in -t_{n}\Omega= -\Omega$ which is a contradiction for we have assumed that $A \cap -\Omega=\emptyset$. This proves our claim. 

Choose $t_0 \geq 1$ such that $t \geq t_0$ implies that $K \cap (A+ta)=\emptyset$. Now note that for $t \geq t_0$, \[||V_{ta}^{*}\phi||^{2}=\int_{A} |\phi(x+ta)|^{2}dx=0.\] Hence $V_{ta}^{*}\phi \to 0$ as $t \to \infty$. This completes the proof. \hfill $\Box$

\begin{ppsn}
\label{CCR flow}
Let $V: P \to B(K)$ be an isometric representation. Then there exists a unique $\en$-semigroup $\alpha^{V}:=\{\alpha^V_{x}\}_{x \in P}$ on $B(\Gamma(K))$ satisfying \[
\alpha^V_{x}(W(u))=W(V_{x}u)~~ \forall ~ x \in P,~u \in K.\]
\end{ppsn}
\textit{Proof.}  We can directly define \begin{equation}\label{CCR defn} \alpha^V_x(T)  = 1_{\Gamma(\ker(V_x^*))} \otimes \Gamma(V_x)T\Gamma(V_x)^* ~~~ \forall T \in B(\Gamma(K)),\end{equation}  where $\Gamma(V_x)$ is the second quantization of $V_x: K \mapsto V_xK$ considered as a unitary operator.  The fact that $V$ is an isometric representation and the linear span of $\{W(u): u \in K\}$ is $\sigma$-weakly dense in $B(\Gamma(K))$ implies that $\{\alpha_{x}\}_{x \in P}$ is a semigroup of endomorphisms. Clearly $\alpha_{x}$ is unital for every $x \in P$. All it remains to check is the continuity property of $\alpha^{V}$.  

Let $\mathcal{A}$ be the linear span of $\{W(u): u \in K\}$. Then $\cla$ is a unital $*$-subalgebra of $B(\Gamma(K))$ and $\mathcal{A}$ is strongly dense in $B(\Gamma(K))$. 
It is clear that for $u,v,w \in K$, the map $P \ni x \to \langle\alpha_{x}(W(u))e(v)|e(w)\rangle$ is continuous. From this, using the fact that exponential vectors are total in $\Gamma(K)$, it is easy to deduce that for $u \in K$ and $\xi,\eta \in \Gamma(K)$, the map $P \ni x \to \langle\alpha_{x}(W(u))\xi|\eta\rangle$ is continuous. Thus it follows that for $A \in \mathcal{A}$, the map $P \ni x \to \langle\alpha_{x}(A)\xi|\eta\rangle$ is continuous for every $\xi,\eta \in \Gamma(K)$. 

Now let $A \in B(\Gamma(K))$ and $\xi,\eta \in \Gamma(K)$ be given. By Kaplansky density theorem, there exists a sequence $A_{n} \in \cla$ such that $A_{n} \to A$ strongly and $||A_{n}|| \leq ||A||$. Since $(A_n)$ is norm bounded, it follows that $A_n \to A$ in the $\sigma$-weak topology. Hence for every $x \in P$, $\alpha_{x}(A_n) \to \alpha_{x}(A)$ in the $\sigma$-weak topology. Thus the sequence of continuous functions $\langle\alpha_{x}(A_n)\xi|\eta\rangle$ converges pointwise to $\langle\alpha_{x}(A)\xi|\eta\rangle$.  This implies that the map $P \ni x \to \langle\alpha_{x}(A)\xi|\eta\rangle$ is measurable. Now Corollary 4.3 of \cite{Murugan_Sundar} implies that for every $A \in B(\Gamma(K))$ and $\xi,\eta \in \Gamma(K)$, the map $P \ni x \to \langle\alpha_{x}(A)\xi|\eta\rangle$ is continuous. This completes the proof. \hfill $\Box$

We call the $\en$-semigroup $\alpha^V$ constructed in Prop.\ref{CCR flow} as the CCR flow associated to the isometric representation $V$. In this paper, we refer to the traditional CCR flow considered in \cite{Arveson} as the $1$-parameter CCR flow.

We end this section by considering tensor product of two $\en$-semigroups. Let $\alpha:=\{\alpha_{x}\}_{x \in P}$ and $\beta:=\{\beta_{x}\}_{x \in P}$ be $\en$-semigroups on $B(H)$ and $B(K)$ respectively. For $x \in P$, there exists a unique normal $*$-endomorphism, denoted $\alpha_{x} \otimes \beta_{x}$, on $B(H \otimes K)$ such that $\alpha_{x}\otimes \beta_{x}(A \otimes B)=\alpha_{x}(A)\otimes \beta_{x}(B)$ for $A \in B(H)$ and $B \in B(K)$.  For the existence of the endomorphism $\alpha_{x} \otimes \beta_{x}$, we refer the reader to Page 21, Paragraph 3 of \cite{Arveson}. It is clear that $\alpha\otimes \beta:=\{\alpha_{x} \otimes \beta_{x}\}_{x \in P}$ is a semigroup of unital normal $*$-endomorphisms of $B(H \otimes K)$. To check the continuity property, one again appeals to a Kaplansky density type argument employed in the proof of Proposition \ref{CCR flow}. 
\begin{rmrk}
The CCR flows admits the following factorisation property. If $V$ and $W$ are isometric representations of $P$ then $\alpha^{V\oplus W}$ is conjugate to $\alpha^{V}\otimes \alpha^{ W}.$
\end{rmrk}

\section{Units} Since  $H^{2}(\bbr_+,\mathbb{T})$ is  trivial, units for $1$-parameter $\en$-semigroups are defined to be semigroups. But we need to bring in the multiplier group while defining units for general $\en$-semigroups over $P$.





\begin{dfn}
Let  $\alpha:=\{\alpha_{x}\}_{x \in P}$ be an $\en$-semigroup on $B(H)$ and $\omega$ be a multiplier on $P$. A strongly continuous family $\{u_{x}\}_{x \in P}$ of bounded operators on $H$ is called a $\omega$-unit if 
\begin{enumerate}
\item[(1)] for $x \in P$ and $T \in B(H)$, $\alpha_{x}(T)u_{x}=u_{x}T$,
\item[(2)] for $x,y \in P$, $u_{x+y}=\omega(x,y)u_{x}u_{y}$, and
\item[(3)] there exists $x \in P$ such that $u_{x} \neq 0$.
\end{enumerate}
If $\omega=1$, we simply call a $\omega$-unit a unit.
\end{dfn}

We denote the collection all $\omega-$units of an  $\en$-semigroup $\alpha$ by $\mathfrak{U}_{\alpha}^{\omega}$ and when $\omega=1$ we simply denote by $\mathfrak{U}_{\alpha}$.

Fix an $\en$-semigroup $\alpha:=\{\alpha_{x}\}_{x \in P}$ on $B(H)$.  Let $\{u_{x}\}_{x \in P}$ be a $\omega$-unit where $\omega$ is a multiplier on $P$. Observe the following.
\begin{enumerate}
\item[(1)] For $x \in P$, $u_{x}^{*}u_{x}$ commutes with every element of $B(H)$. Thus for every $x \in P$, $u_{x}^{*}u_{x}$ is a scalar.
\item[(2)] Let $x \in P$ be such that $u_{x}\neq 0$. Then $u_{x}u_{0}=\overline{\omega(x,0)}u_{x}$. Premultiplying by $u_{x}^{*}$, we obtain that $u_{0}=\omega(x,0)$. Thus $u_0$ is a non-zero scalar multiple of identity.
\item[(3)] Let $\{U_{x}\}_{x \in P}$ be an $\alpha$-cocycle and let $\beta:=\{Ad(U_x)\circ \alpha_{x}\}_{x \in P}$. Then $\{U_{x}u_{x}\}_{x \in P}$ is a $\omega$-unit for $\beta$.
\end{enumerate}

\begin{ppsn}\label{equivalent-units}
Let $\alpha:=\{\alpha_{x}\}_{x \in P}$ be an $\en$-semigroup on $B(H)$. Let $\omega_1$ and $\omega_2$ be multipliers on $P$. Suppose that $\alpha$ admits a $\omega_1$-unit and also a $\omega_2$-unit then 
$[\omega_1]=[\omega_2]$.
\end{ppsn}

\textit{Proof.} Let $\{u_{x}\}_{x \in P}$ be a $\omega_1$ unit  and $\{v_{x}\}_{x \in P}$ be a $\omega_2$ unit for $\alpha$. For  $x \in P$, let  
$\bbc\ni f(x) =v_{x}^{*}u_{x}$, then, by $(2)$ of above observations,  $f(0)$ is  non-zero. The continuity of $f$ follows from the continuity of $\{u_x\}_{x \in P}$ and $\{v_{x}\}_{x \in P}$. Now for $x,y \in P$,
\begin{align*}
f(x+y)&=v_{x+y}^{*}u_{x+y} \\
         &=\overline{\omega_{2}(x,y)}\omega_{1}(x,y)v_{y}^{*}v_{x}^{*}u_{x}u_{y} \\
         &=\overline{\omega_{2}(x,y)}\omega_{1}(x,y)f(x)v_{y}^{*}v_{y} \\
          &=\overline{\omega_{2}(x,y)}\omega_{1}(x,y)f(x)f(y) 
\end{align*}
The above calculation, the fact that $P$ is a closed convex cone, $f$ is continuous and $f(0) \neq 0$ implies that $f(x) \neq 0$ for every $x \in P$. Also  $f(x)f(y)f(x+y)^{-1} \in \mathbb{T}$, for all $x , y \in P$. For $x \in P$, set $g(x)=\frac{f(x)}{|f(x)|}$. Then 
the preceding calculation shows that $\omega_{1}(x,y)g(x)g(y)g(x+y)^{-1}=\omega_{2}(x,y)$ for all $x,y\in P$. This implies that $[\omega_1]=[\omega_2]$. This completes the proof. \hfill $\Box$

Let $A \in B_{0}(d)$ be given. We denote the map $\mathbb{R}^{d} \times \mathbb{R}^{d} \ni (x,y) \to e^{i\langle Ax|y\rangle} \in \mathbb{T}$ also by $\omega_{A}$. For $x \in P$, let $U_{x}^{A}$ be the unitary operator on $L^{2}(\mathbb{R}^{d})$ defined by 
\begin{equation}
U_{x}^{A}(f)(y):= \omega_{A}(x,y-x)f(y-x)  
\end{equation}
for $f \in L^{2}(\mathbb{R}^{n})$. It is clear that $\{U_{x}^{A}\}_{x \in P}$ is a strongly continuous family of unitaries such that for $x,y \in P$, $U_{x}^{A}U_{y}^{A}=\omega_{A}(x,y)U_{x+y}^{A}$. For $x \in P$, let $\alpha_{x}^{A}:=Ad(U_{x}^{A})$. Then $\alpha^{A}:=\{\alpha_{x}^{A}\}_{x \in P}$ is an $\en$-semigroup on $B(L^{2}(\mathbb{R}^{d}))$. Note that $\alpha_{x}^{A}$ is an automorphism for each $x \in P$ and $\alpha^A$ can be extended to an automorphism group on $\mathbb{R}^d$.

Let $\mathfrak{E}$ denote the set of equivalence classes of $\en$-semigroups on $B(L^{2}(\mathbb{R}^{d}))$ which can be extended to automorphism groups. Here the equivalence relation is that of cocycle conjugacy. The following Proposition can be concluded from  Proposition \ref{equivalent-units}, Wigner's Theorem, discussions in Section 1.4 and Theorem 3.6.6 in  \cite{KRP1}.

\begin{ppsn}
\label{classification}
With the foregoing notations, the map $B_{0}(d) \ni A \to [\alpha^{A}] \in \mathfrak{E}$ is a bijection. 
\end{ppsn}

\begin{rmrk}
Let $\alpha:=\{\alpha_{x}\}_{x \in P}$ be an $\en$-semigroup on $B(H)$. Let $\omega_1$ and $\omega_2$ be multipliers on $P$ which are cohomologous . Suppose that $\alpha$ admits a $\omega_1$-unit. Let $f:P \to \mathbb{T}$ be a continuous function such that $\omega_{2}(x,y)=f(x+y)f(x)^{-1}f(y)^{-1}\omega_{1}(x,y)$.  
Then it is easily verifiable that the map 
\[
\mathfrak{U}_{\alpha}^{\omega_1} \ni \{u_{x}\}_{x \in P} \to \{f(x)u_{x}\}_{x \in P} \in \mathfrak{U}_{\alpha}^{\omega_2}
\]
is a bijection.
\end{rmrk}

\begin{xmpl} Let $V$ be an isometric representation of $P$ and $\alpha^V$ be the associated CCR flow.  It can be directly verified from equation \ref{CCR defn} that $\{\Gamma(V_{x})\}_{x \in P}$ is a unit for $\alpha^{V}$. 

Let $A \in B_{0}(d)$ be given. Denote the multiplier $P \times P \ni (x,y) \to e^{i\langle Ax|y\rangle}$ by $\omega_{A}$. Let  $\{U_{x}^{A}\}_{x \in P}$ be the unitaries  and $\alpha^{A}$ be the $\en$-semigroup considered in Prop.\ref{classification}.  Then $\{\Gamma(V_{x}) \otimes U_{x}^{A}\}_{x \in P}$ is a $\omega_{A}$-unit for the $\en$-semigroup $\alpha^{V} \otimes \alpha^{A}$. But $\alpha^{V} \otimes \alpha^{A}$ is not a pure $\en$-semigroup. The problem of producing a pure $\en$-semigroup admitting a non-trival $\omega_{A}-$unit  is open.
\end{xmpl}

The set of all units of a  CCR flow $\alpha^V$ is described by the set of additive cocycles of $V$, which are defined  as below.  These are supposed to be called more precisely as local additive cocycles, but since we use  only local additive cocycles, we just call them as additive cocycles.

\begin{dfn}
Let $V: P \to B(K)$ be a strictly continuous isometric representation.
A continuous function $h : P \to K$ is called a $P$-additive cocycle  for $V$  if 
\begin{enumerate}
 \item[(i)] for all $ x \in P$,  $h_x \in \ker(V_x^*)$, and
 \item[(ii)] for all $ x, y \in P$,  $h_x+ V_x h_y = h_{x+y}.$
\end{enumerate} 

An $\Omega$-additive cocycle is defined by replacing $P$ by $\Omega$ in above definition. For a $P$-additive cocycle $\xi$,  the restriction to $\Omega$ is clearly an $\Omega$-additive cocycle.  
\end{dfn}

We denote by $\mathscr{A}^V(P)$ (respectively $\mathscr{A}^V(\Omega))$ the set of all $P-$additive cocycles (respectively $\Omega-$additive cocycles ) of $V$.  Clearly $\mathscr{A}^V(P)$ forms a vector space with respect to natural addition and scalar multiplication.  We endow $\mathscr{A}^V(P) $ with the topology of uniform convergence in norm on all compact subsets of $P$.

The following Lemma is well known to those who has  had a  faint exposure to  Ore semigroups. We include the proof for completeness.

\begin{lmma} \label{linext}
\noindent  Let $G$ be a topological group and $\phi: \Omega \to G$ be a continuous map such that $\phi(a+b)=\phi(a)\phi(b)$ for every $a,b \in \Omega$. Then $\phi$ extends to a unique continuous group homomorphism from $\mathbb{R}^{d}$ into $G$.
\end{lmma}
\textit{Proof.} Let $x \in \mathbb{R}^{d}$ be given. Write $x=a-b$ with $a,b \in \Omega$. Define $\widetilde{\phi}(x)=\phi(a)\phi(b)^{-1}$. We claim that $\widetilde{\phi}$ is well defined. Suppose $a_1-b_1=a_2-b_2$ where $a_1,a_2,b_1,b_2 \in \Omega$. Then $a_1+b_2=a_2+b_1$. Hence $\phi(b_2)\phi(a_1)=\phi(a_2)\phi(b_1)$. Hence $\phi(a_1)\phi(b_1)^{-1}=\phi(b_{2})^{-1}\phi(a_2)=\phi(a_{2})\phi(b_{2})^{-1}$. This shows that $\widetilde{\phi}$ is well defined. Clearly $\widetilde{\phi}$ is continuous on $\Omega-b$ for every $b \in \Omega$. Since $\{\Omega-b: b \in \Omega\}$ is an open cover of $\mathbb{R}^{d}$, it follows that $\widetilde{\phi}$ is continuous. This completes the proof. \hfill $\Box$

\begin{lmma}\label{boundedh}
Let $h:\Omega \to K$ be an $\Omega$-additive cocycle. Then there exists $\mu\in \bbr^d$ such that $||h_a||^2=\langle \mu | a\rangle $ for every $a\in \Omega.$
 \end{lmma}
\textit{Proof}. Define a function $f:\Omega \to \bbr$ by the formula:
\[ f(a)=\langle h_a| h_a\rangle,\;\; \text{for every } a\in \Omega.\]

Let $a, b \in \Omega $ be given. Then 
\begin{align}
 f(a+b)&=\langle \xi_{a+b}| \xi_{a+b} \rangle \cr
 &= \langle \xi_a + V_a \xi_b | \xi_a + V_a \xi_b \rangle \cr
 &=\langle \xi_a | \xi_a \rangle + \langle V_a \xi_b | V_a \xi_b \rangle  \;\quad \left(\text{Since } \langle \xi_a | V_a \xi_b \rangle = 0.\right)\cr
 &=\langle \xi_a | \xi_a \rangle + \langle \xi_b | \xi_b \rangle = f(a)+f(b).\cr
\end{align}
Continuity of $f$ follows from continuity of $h$.
Then by lemma \ref{linext}, there exists $\mu \in \bbr^d$  such that for every $a \in \Omega$, $f(a)=\langle a |\mu  \rangle.$ This implies $||h_{a}||^2=\langle a |\mu \rangle.$

\begin{lmma} \label{Paddcocycle}
Let $V:P \to B(K)$ be a strictly continuous isometric representation. If $h :\Omega \to K$ is a $\Omega$-additive cocycle, then there exists an unique extension $\tilde{h}: P \to K$ of $h$ such that $\tilde{h}$ is a $P$-additive cocycle.
 \end{lmma}
 
\textit{Proof.} Let $U_a=W(h_a)$ for $ a \in \Omega$. Then $U_a$ satisfies the gauge cocycle conditions on $\Omega$. By Prop. 4.4 of  \cite{Murugan_Sundar}, 
there exists an extension $\{\widetilde{U_x}\}_{x \in P}$ of $U:=\{U_a\}_{a \in \Omega}$, as an $\alpha$-cocycle. The strong continuity implies  that $\{\widetilde{U_x} \}_{x \in P}$ is a gauge cocycle. 

Let $x \in P$. Choose a sequence $\{a_n\}$ in $\Omega$ such that $a_n$ converges to $x$. We claim that  the sequence $\{h_{a_n}\}$ converges and the limit  is independent of the chosen sequence. 
Since $a_n$ converges to $x$, by strong continuity, $W(h_{a_n})e(0)=e^{-\frac{1}{2} \|h_{a_n}\|^2}(e(h_{a_n}))$ converges to $\widetilde{U_x}e(0)$.  By projecting onto the $0-$particle space and onto the $1-$particle space, we find that $e^{-\frac{1}{2} \|h_{a_n}\|^2}$ and $e^{-\frac{1}{2} \|h_{a_n}\|^2}h_{a_n}$ converges. By Lemma \ref{boundedh} the sequence $\{h_{a_n}\}$ is bounded and  
 $h_{a_n}$ converges.
Define $\widetilde{h_x} := \lim_{n \to \infty} h_{a_n}.$   Suppose $b_n \in \Omega$ be another sequence such that $b_n$ also converges to $x$ and $\widetilde{h'_x}:= \lim_{n \to \infty} h_{b_n}.$  Let $ c \in \Omega$, then both $\Omega \ni a_n +c, b_n +c $ converge to $x +c \in \Omega $. This implies both $h_{a_n +c}$ and  $h_{b_n+c}$ converge to $h_{x+c}.$ But $h_{a_n +c}= h_{a_n} + V_{a_n}h_c \to \widetilde{h_x} + V_xh_c$ and  $h_{b_n +c}= h_{b_n} + V_{b_n}h_c \to \widetilde{h'_x} + V_xh_c$. Hence $\widetilde{h_x}=\widetilde{h'_x}.$   By the same technique we can see $x \to h_x$ is continuous and it is also easy to check that $\widetilde{h}$ is a $P$-additive cocycle.
\hfill $\Box$

 For a strictly continuous isometric representation $V: P \to B(K)$ and a  $\xi_x \in \Gamma(\ker(V_x^*))$ with $x \in P$,
 define $R^{\xi_x}$  by $R^{\xi_x} (\eta) =\xi_x \otimes \Gamma(V_x)\eta$  for all $\eta  \in \Gamma(K).$  Then  it is easy to check that $R^{\xi_x}  \in   B(\Gamma(K))$  and 
 $\alpha_{x}(X)R^{\xi_x}=R^{\xi_x} X ~~ \forall X \in B(\Gamma(K)).$

For $\mu \in \bbc^n$ and $h \in \mathscr{A}^V(P)$,  define $T^{\mu, h}= \{T^{\mu, h}_x\}_{x \in P}$ with $T^{\mu, h}_x = e^{\langle x, \mu\rangle} R^{e(h_x)}.$ 
 Since $\Gamma(V_x)$ is a 1-unit for $\alpha^V$, thanks to Proposition \ref{equivalent-units}, there does not exist any $\omega-$unit for $\alpha^V$ with $\omega$ not equivalent to $1$.  The units of the CCR flow $\alpha^V$ admits the following description.

\begin{thm} \label{CCR units}
 Let $\alpha^V$  be the CCR flow  associated to a pure  strictly continuous semigroup of isometric representation $V$. Then the map 
 $\bbc^d \times \mathscr{A}^V(P) \ni (\mu,h) \to T^{\mu,h} \in  \mathfrak{U}_{\alpha}$ is a bijection.
 \end{thm}	
 
  \textit{Proof.}   It is easy to verify that $T^{\mu, h}$ is a unit for  $\alpha^V$, for any $(\mu, h) \in \bbc^n \times \mathscr{A}^V(P)$.  It is also a direct verification to see that the above map is   injective, by evaluating on vacuum vectors.
  
  To prove surjectivity, let $\{u_x\}_{x\in P}$ be a unit for $\alpha^V$.
  For any fixed $a \in \Omega$, the  $\en$-semigroup $\{\alpha^V_{ta}: t \in \bbr_+ \}$ is the 1-parameter CCR flow 
   associated with the pure family of isometries $\{V_{ta}: t \in \bbr_+\}$. Hence, thanks to  \cite[Theorem 2.6.4]{Arveson},  there exists $\lambda(a) \in \bbc^{*}$ and $h(a) \in \ker(V_a^*)$ such that 
   $u_a  =\lambda(a) R^{e(h_a)}.$ 
   For any $a,b \in \Omega$, the relation $u_au_b = u_{a+b}$ evaluated on  vacuum vectors implies that 
 $$\lambda(a+b) = \lambda(a)\lambda(b), ~~ h_{a+b}= h_a + V_ah_b.$$
   The continuity of $\lambda$ and $\{h_{a}\}_{a \in \Omega}$ follows from the continuity of $\{u_{a}\}_{a \in \Omega}$. Now the rest of the proof of surjectivity follows from Lemmas \ref{linext} and Lemma \ref{Paddcocycle}.
  \hfill $\Box$
  
    A proof of the following fact is contained in  \cite[Theorem 2.6.4]{Arveson}.
  
  \begin{ppsn}
  Additive cocycles of one parameter  right shift $\{S_t\}_{t \in \bbr_+}$ on $L^2(\bbr_+, \k)$ are given by $h_t = 1_{(0,t)} \otimes k$ for some $k \in \k$. 
  \end{ppsn}
  
It is an elementary exercise to prove the first statement of the following Lemma, from which the second statement follows easily.

\begin{lmma}\label{3tensors}  (i) Let $K_1, K_2, K_3$ be complex Hilbert spaces  and $U_\sigma$ be the unitary operator which flips the first and second component in the tensor product $K_1 \otimes K_2\otimes K_3$. If 
 $f_1 \otimes F_{23}  = U_\sigma (g_2\otimes G_{13})$ for some $f_1 \in K_1, F_{23} \in K_2\otimes K_3$ and $g_2 \in K_2,G_{13} \in K_1\otimes K_3$, then  there exists a $f_3 \in K_3$ such that $$f_1 \otimes F_{23} =f_1 \otimes g_2 \otimes f_3 =  U_\sigma (g_2\otimes G_{13}).$$

\noindent (ii) Suppose $R_1 \otimes R_{23} = U_\sigma (R_2\otimes R_{13})U^*_\sigma$ for $R_1 \in B(K_1), R_{23} \in B( K_2\otimes K_3)$ and $R_2 \in B(K_2), R_{13} \in B( K_1\otimes K_3)$, then there exits an $R_3 \in B(K_3)$ such that
$$R_1 \otimes R_{23} = R_1 \otimes R_2 \otimes R_3 =U_\sigma (R_2\otimes R_{13})U^*_\sigma.$$
\end{lmma}

\begin{ppsn}\label{0addits}
  Let $P = \bbr_+^d$ with $d \geq 2$ and let $A = \times_{i=1}^d \bbr_i$ where $\bbr_i \in \{\bbr, \bbr_+\}$ with at least one $\bbr_i = \bbr_+$. Then the semigroup of shifts $\{S^A_x\}_{x \in P}$ on $L^2(A, \k)$ does not admit any non-trivial additive cocycles.  
  
  Also the semigroup of shifts $\{S^{A_T^{(a,b)}}_x\}_{x \in \bbr_+^2}$, discussed in Example \ref{A12}, does not admit  any non-trivial additive cocycles.
  \end{ppsn}
  
 \textit{Proof.} Suppose $\bbr_{i_0} =\bbr_+$ for some $i_0$, in the direct product  $A = \times_{i=1}^d \bbr_i$,  and 
 $h=\{h_{x}\}_{x \in P}$ be an additive cocycle for the semigroup of shifts $\{S^A_x\}_{x \in P}$ on $L^2(A, \k)$. 
  Let $e_{i_0}=(0,...,1,..0),$ with 1 in the $i_0$-th place and $0$ elsewhere.  Then for any $s \in \bbr_+$, we have $h_{se_{i_0}}=U_{\sigma_{i_0}}\left(1_{(0,s)}\otimes f_{i_0}\right)$ for some $f_{i_0} \in K({i_0})$, where  $U_{\sigma_i}$ is the unitary which flips ${i_0}-$th tensor component with the first component in $ \otimes_{i=1}^d L^2(\bbr_i) \otimes \k$, and $K({i_0}) = \otimes_{i=1}^{{i_0}-1}L^2(\bbr_i) \otimes \widehat{L^2(\bbr_+)} \otimes_{i={i_0}+1}^d L^2(\bbr_i)\otimes \k$. Here $\widehat{.}$ indicates that the ${i_0}$-th component Hilbert space in the tensor product is removed. 
 
 Suppose for $\bbr_j = \bbr$, notice $h_{te_j} =0$ for all $t \in \bbr_+$. Then for  $s,t \in \bbr_+$, the relation $h_{se_{i_0}} + S^A_{se_{i_0}} h_{te_j} = h_{se_{i_0}+te_j}= h_{te_j} + S^A_{te_j} h_{se_{i_0}}$,   implies that $U_{\sigma_{i_0}}\left(1_{(0,s)}\otimes f_{{i_0}}\right) = S^A_{te_j} \left(U_{\sigma_{i_0}}\left(1_{(0,s)}\otimes f_{{i_0}}\right)\right)$ which consequently imply that  $f_{i_0}=0$.
 On the other hand, suppose $\bbr_j = \bbr_+$ for some $j\neq {i_0}$, then $h_{te_j}=U_{\sigma_j}\left(1_{(0,t)}\otimes f_{j}\right)$ for some $f_j \in K(j)$. Set $E_x = 1-S^A_x(S^A_x)^*$ for any $x \in P$.  Notice $E_{se_{i_0}}$ and $E_{t e_j}$ commute for any  $s, t \in \bbr_+$.
 Now  applying $E_{se_{i_0}}E_{te_j}$ on both sides of $h_{se_{i_0}} + S^A_{se_{i_0}} h_{te_j} = h_{te_j} + S^A_{te_j} h_{se_{i_0}}$, we get  $$E_{te_j}U_{\sigma_{i_0}}\left(1_{(0,s)}\otimes f_{{i_0}}\right) = E_{se_{i_0}} U_{\sigma_j}\left(1_{(0,t)}\otimes f_{j}\right)= U_{\sigma_j}\left(1_{(0,t)}\otimes E_{se_{i_0}} f_{j}\right)~~ \forall s, t \in \bbr_+.$$  (Here and in the rest of the paragraph we have slightly abused the notation.) When taking limit $t\rightarrow \infty$, the left hand side goes to $U_{\sigma_{i_0}}\left(1_{(0,s)}\otimes f_{{i_0}}\right)$. But as $t\rightarrow \infty$, the right hand side goes to $U_{\sigma_j}\left(1_{(0,\infty)}\otimes E_{se_{i_0}} f_{j}\right)$ which is not square integrable, unless $E_{se_{i_0}}f_j =0$. Since this holds for all $s$, 
 we have  $f_j =0$ and hence also   $f_{i_0} =0$.  Since ${i_0},j$ were arbitrary we conclude that $h_{x}=0$ for all $x \in P$. The proof of the first part of the Lemma is over.

 Now we consider the case of $\{S^{A_T^{(a,b)}}_x\}_{x \in R_+^2}$ of Example \ref{A12}. Since $(a,b)\in \bbr^2$ is fixed we simply write $\{S^{A_T^{(a,b)}}_x\}_{x \in R_+^2}$ by $\{S^{A_T}_x\}_{x \in R_+^2}$. 
We set  $A_{T1}= [0,\infty) \times [0,b)$ and $A_{T2}= [a, \infty) \times [b, \infty)$, so that $A_T = A_{T1}  ~\bigsqcup ~A_{T2}$. We also set  $A_{1T} = [0,\infty) \times [0,\infty)$  and $A_{2T} = [a,0) \times [b, \infty)$ so that $A_T = A_{1T}  ~\bigsqcup ~A_{2T}$.
Let $E_1$ and $E_2$ be the projections onto $L^2(A_{T1})$  and $L^2(A_{T2})$  respectively, and let $F_1$ and $F_2$ be the projections onto $L^2(A_{1T})$  and $L^2(A_{2T})$  respectively.  We denote  the one parameter shifts on first component of  $L^2( (0,\infty))\otimes L^2((0,b))\otimes k$ and $L^2( (a,\infty))\otimes L^2((b, \infty))\otimes k$ by $S^{1}_t$ and $S^{2}_t$ respectively. Also we denote the one parameter shifts on the second component of $L^2( (0,\infty))\otimes L^2((0, \infty))\otimes k$ and $L^2( (a,0))\otimes L^2((b, \infty))\otimes k$ by $S^{3}_t$ and $S^{4}_t$respectively.

Suppose $h =\{h_{(s,t)}:s, t, \in \bbr_+\}$ is an additive cocycle for 
$S_{A_T}$, then both $\{E_1 h_{(s,0)}\}$ and  $\{E_2 h_{(s,0)}\}$
form additive cocyles for $S^{1}_t$ and $S^{2}_t$ respectively.  Similarly, $\{F_1 h_{(0,t)}\}$ and  $\{F_2 h_{(0,t)}\}$
form additive cocyles for   $S^{3}_t$ and $S^{4}_t$respectively. Hence there exists $ f_1 \in L^2(0,b)\otimes \k ,~f_2 \in L^2(b,\infty)\otimes \k,~f_3 \in L^2(a,0)\otimes \k,\; f_4 \in L^2(0,\infty)\otimes \k$ satisfying
 \[h_{(s,0)}=1_{(0,s)}\otimes f_1 +1_{(a,a+s)}\otimes f_2;~~~
  h_{(0,t)}=f_3 \otimes 1_{(b,b+t)}+f_4 \otimes 1_{(0,t))}.\]
 (Be cautioned that the subscript in $h_{(\cdot, \cdot)}$ is a tuple in $\bbr_+^2$ and the subscript in $1_{(\cdot, \cdot)}$ is an interval. Also notice we are slightly abusing notations, by not writing a tensor flip explicitly in $f_3 \otimes 1_{(b,b+t)}+f_4 \otimes 1_{(0,t))}$. )We have 
  \begin{align}\label{vectorlhs}
 h_{(s,t)}&=h_{(s,0)}+ S^{A_T}_{(s,0)}h_{(0,t)}\cr
 &=1_{(0,s)}\otimes f_1 +1_{(a,a+s)}\otimes f_2  +  S^{A_T}_{(s,0)} \left(f_3 \otimes 1_{(b,b+t)} +f_4 \otimes 1_{(0,t)}\right)\cr
&=1_{(0,s)}\otimes f_1 +1_{(a,a+s)}\otimes f_2 +S_s^{1} E_1(f_4 \otimes 1_{(0,t)})+ S_s^{2}E_2(f_4 \otimes 1_{(0,t)})    + S_s^{2}(f_3 \otimes 1_{(b,b+t)}).\cr
                                                      \end{align}
                                                      On the other hand                                                        
\begin{align} \label{vectorrhs}                                      
h_{(s,t)}&=h_{(0,t)}+ S^{A_T}_{(0,t)} h_{(s,0)}\cr
 &=f_3\otimes 1_{(b,b+t)} +f_4\otimes 1_{(0,t)} + S^{A_T}_{(0,t)} (1_{(0,s)}\otimes f_1 +1_{(a,a+s)}\otimes f_2)\cr
&=f_3\otimes 1_{(b,b+t)} +f_4\otimes 1_{(0,t)} +  S_t^4(1_{(0,s)}\otimes f_1)
+S_t^3 F_1(1_{(a,a+s)}\otimes f_2)+S_t^4 F_2(1_{(a,a+s)}\otimes f_2).\cr
                                                     \end{align}

Fix $t= b$. Let $P_{[0,s)}\otimes I_{(0,b)} \otimes I$ be the projection of $L^2(0,\infty)\otimes L^2(0,b) \otimes \k$ onto $L^2(0,s)\otimes L^2(0,b) \otimes \k$. 
        Now
applying $\left(P_{[0,s)}\otimes I_{(0,b)} \otimes I\right) E_1$ in equations \ref{vectorlhs} and \ref{vectorrhs} , we get
\[1_{(0,s)}\otimes f_1= \left(P_{[0,s)}\otimes I_{(0,b)}  \right)f_4 \otimes 1_{(0,b)}.\] Thanks to Lemma \ref{3tensors}, \[ \left(P_{[0,s)}\otimes I \right)f_4 =1_{(0,s)}\otimes k,\;\;\text{and} \;f_1=1_{(0,b)}\otimes k,\;\;\text{for some}\;\; k\in \k, \;\; \forall\; s\in \bbr_+.\]
By letting $s\to \infty,$ we get $f_4=1_{\bbr_+}\otimes k$, but since $f_4$ is square integrable $k=0$. Therefore $f_4=0$ and $f_1=0$.

Now we have from Equations \ref{vectorrhs} and \ref{vectorlhs}
$$1_{(a,0)}\otimes f_2  + S_s^{2}(f_3 \otimes 1_{(b,b+t)}) =f_3\otimes 1_{(b,b+t)}  + S^{A_T}_{(0,t)}(1_{(a,a+s)}\otimes f_2).$$ 
Let $I_{(a,0)}\otimes P_{[b,b+t)} \otimes I$ be the projection of $L^2(a,0)\otimes L^2(b, \infty) \otimes \k$ onto $L^2(a,0)\otimes L^2(b, b+t) \otimes \k$.  Now applying $\left(I_{(a,0)}\otimes P_{[b,b+t)} \otimes I\right) F_2$ on both sides, for $s \geq -a$, we get  $$1_{(a,0)} \otimes (P_{(b, b+t)}\otimes I) f_2 = f_3 \otimes 1_{(b,b+t)}~~ \forall t \in \bbr_+ .$$ By the same  argument we get $f_2=0, f_3=0,$ and therefore $h_{(s,t)}=0,\;\forall \; s,t\in\bbr_+.$       \hfill $\Box$ 

\medskip

The following theorem, which is immediate from the above discussions, asserts that all our examples arising from CCR are of type II$_0$, that is they admit only one unit up to scalars.

\begin{thm}
   Let $P = \bbr_+^d$ and let $A = \times_{i=1}^d \bbr_i$ where $\bbr_i \in \{\bbr, \bbr_+\}$ with at least one $\bbr_i = \bbr_+$.  Then for the $\en$-semigroup $\alpha^{S^A}$ associated with the isometric representation $S^A$, $$\mathfrak{U}(  \alpha^{S^A}) = \{ \{e^{\langle x, \mu\rangle} \Gamma(S_x)\}_{x \in P}: \mu \in \bbc^n\}.$$
  Also for the $\en$-semigroup $\{\alpha^{S^{A_T}}\}_{x \in R_+^2}$ (discussed in Example \ref{A12}),
  $$\mathfrak{U}(\alpha^{S^{A_T}}) = \{ \{e^{\langle x, \mu\rangle} \Gamma(S^{A_T})\}_{x \in P}: \mu \in \bbc^n\}.$$ 
 \end{thm}
 
 \begin{rmrk} Let $P = \bbr_+^d$.  Let $\alpha_i$ be $1-$parameter $\en$-semigroups $\alpha_i$ on $B(H_i)$ for $i=1,\cdots d$. Then we can define an  $\en$-semigroup $\alpha = \prod_{i=1}^d \alpha_i$ on $B(\otimes_{i=1}^d H_i)$  by  $$\alpha_{(t_1, t_2\cdots t_d)}(\otimes_{i=1}^d X_i)) = \otimes_{i=1}^d \alpha_{t_i}(X_i), ~~ \forall (t_1, t_2\cdots t_d) \in \bbr_+^d.$$
 We say an $\en$-semigroup $\alpha$  is factorizable if there exists $1-$parameter $\en$-semigroups $\alpha_i$  for $i=1\cdots d$, such that $\alpha$ is cocycle conjugate to $\prod_{i=1}^d\alpha_i$.

For the  CCR flows $\alpha^{S^A}$ and $\alpha^{S^{A_T}}$,  all the restricted $1-$parameter $\en$-semigroups are of type I, whereas the $\en$-semigroup over $P$  is of type II$_0$. Therefore it follows that the $\alpha^{S^A}$ and $\alpha^{S^{A_T}}$ are not factorizable.  (Readers may refer to \cite{Arveson} for the definitions of type I and II, and also it is an easy verification to check  tensor product of type I  $\en$-semigroups is type I.)
 \end{rmrk}
 
 \section{Standard form and conjugacy} 
\begin{dfn}
An $\en$-semigroup $\alpha$ on $B(H)$ is said to be in the standard form if there exists  a $\Omega \in H$ such that $$\langle \Omega | \alpha_x (X) \Omega\rangle =  \langle \Omega |  X \Omega\rangle ~~\forall x\in P,~\forall X \in B(H).$$
That is there exists a pure invariant state for $\alpha$.
\end{dfn}
 
When  $\alpha$ is in the standard form, there always exists a canonical unit $\{T_x\}_{x \in P}$ for $\alpha$, defined  
as follows. When $X\Omega =Y\Omega$ for $X, Y \in B(X)$, $$\|\alpha_x(X)\Omega -\alpha_x(Y) \Omega\|^2 = \langle \Omega, \alpha_x(X^*X - X^*Y-Y^*X +Y^*Y)\Omega\rangle =  \|X\Omega -Y \Omega\|^2=0.$$ Now define  $T_x X \Omega = \alpha_x(X) \Omega,$ for all $x \in P$. Then $T_x$ is a well-defined isometry and $T =\{T_x: x \in P\}$ is a unit for $\alpha$. We call $T$ as the canonical unit for $\alpha$.

 The following Theorem generalizes a result, known for $1-$parameter  $\en$-semigroups, to $\en$-semigroups over $P$.
 
 \begin{thm}\label{standard form} Let $\alpha$ and $\beta$ be pure $\en$-semigroups in standard form. We also assume that the gauge group action (by left multiplication) on the set of all units is transitive for $\alpha$. Then $\alpha$ is cocycle conjugate to $\beta$ if and only if $\alpha$ is conjugate to $\beta$.
 \end{thm}
 
  \textit{Proof.} By replacing with a conjugate $\en$-semigroup, if needed, we assume that both $\alpha$ and $\beta$ are acting on same $B(H)$, and that $\beta$ is a cocycle perturbation of $\alpha$. Let $\Omega_\alpha, \Omega_\beta\in H$ be the invariant (unit) vector states for $\alpha$, $\beta$ respectively, and let $T^\alpha, T^\beta$ be the corresponding canonical units. Since, the action of gauge group is transitive on units, we assume that there exists  $\{U_x: x \in P\}$, a unitary cocycle for $\alpha$,   satisfying $$\beta_x = Ad(U_x)\alpha_x; ~~U_x T^\alpha_x = T^\beta_x, ~~\forall x \in P.$$ 
  
  Denote  $\cla_x = \alpha_x(B(H))^\prime, \clb_x =  \beta_x(B(H))^\prime$.  Then $\cla_x \subseteq \cla_y, \clb_x\subseteq \clb_y$ if $x \leq y$, and since both $\alpha$ and $\beta$ are pure, we have $$\bigvee_{x \in P} \cla_x = B(H) = \bigvee_{x \in P} \clb_x .$$ It is easy to verify that $\theta_x= Ad(U_x)|_{\cla_x}$ is an isomorphism between $\cla_x$ and $\clb_x$ for all $x \in P$. Further, for $x \leq y$, and $T \in \cla_x$, $$U_yTU_y^* = U_x \alpha_x(U_{y-x}) T \alpha_x(U_{y-x})^*U_x^*=U_xT\alpha_x\left(U_{y-x}U_{y-x}^*\right)U_x^* = U_xT U_x^*,$$ and hence $\theta_y|_{\cla_x} = \theta_x$. In particular if $T \in \cla_x\cap \cla_y$ for some $x, y \in P$, then $$\theta_x(T) = \theta_{x+y}(T) =\theta_y(T).$$  Therefore, if we define $U T \Omega_\alpha = \theta_x(T)\Omega_\beta$ when $T \in \cla_x$, then $U$ is well-defined. 
  
  Notice for any $T \in \cla_x$, $(T^\alpha_x)^* T T^\alpha_x$ commutes with all operators in $B(H)$, hence is a scalar multiple of  the identity operator. For the same reason, for any $S \in \clb_x$, $(T^\beta_x)^* S T^\beta_x$  is also a scalar multiple of  the identity operator. We denote that both these scalar by $\langle  T T^\alpha_x,  T^\alpha_x \rangle_x$ and $\langle  S T^\beta_x,  T^\beta_x \rangle_x$ respectively. 
  Now for $S \in \cla_x,  T  \in \cla_y$ we have  
  \begin{align*} \langle  \theta_{x}(T)\Omega_\beta |  \theta_{y}(S)\Omega_\beta \rangle   & =  \langle  U_{x+y}S^*T U_{x+y}^*\Omega_\beta | \Omega_\beta \rangle 
   =   \langle  S^*T U_{x+y}^* T_{x+y}^\beta \Omega_\beta | U_{x+y}^*  T_{x+y}^\beta \Omega_\beta  \rangle  \\
  & =  \langle  S^*T T_{x+y}^\alpha \Omega_\beta |  T_{x+y}^\alpha \Omega_\beta  \rangle 
   =   \langle  S^*T T_{x+y}^\alpha |  T_{x+y}^\alpha \rangle_{x+y} \|\Omega_\beta\|^2\\
  & =  \langle  S^*T T_{x+y}^\alpha  |  T_{x+y}^\alpha \rangle_{x+y}  \|\Omega_\alpha\|^2
   =  \langle  S^*T T_{x+y}^\alpha  \Omega_\alpha |  T_{x+y}^\alpha  \Omega_\alpha \rangle 
   =  \langle T  \Omega_\alpha , S  \Omega_\alpha \rangle
 \end{align*} 
  So $U$ preserves inner products and extends as a unitary map on $H$. 
  
  For $T \in \cla_x, S \in \clb_y,$ \begin{align*} Ad(U)(T) S\Omega_\beta & = UTU^*S\Omega_\beta= UTU_y^*SU_y\Omega_\alpha = U_{x+y}TU_y^*SU_yU_{x+y}^*\Omega_\beta\\ 
  &= U_{x+y}TU_y^*SU_y U_y^*\beta_y(U_x^*)\Omega_\beta = U_{x+y}TU_{x+y}^*S \Omega_\beta.
  \end{align*}  So we have $Ad(U)(T) = \theta_{x+y}(T) = \theta_x(T)$ for all $T \in \cla_x$. 
  
 To complete the proof it is enough if we prove that $$Ad(U)\alpha_x Ad(U^*) = Ad(U_x)\alpha_x~~\forall x \in P.$$ Thanks to the purity of $\beta$, it is enough to verify this relation for arbitrary $T \in \clb_y$. Indeed, \begin{align*}Ad(U)(\alpha_x ((Ad(U^*))(T))) & = Ad(U)(\alpha_x ((Ad(U_y^*))(T)))= Ad(U_{x+y})(\alpha_x ((Ad(U_y^*))(T)))\\ &= U_x\alpha_x(U_y) \alpha_x (U_y^*TU_y)\alpha_x(U_y) ^*U_x^*  = U_x \alpha_x(T) U_x^*~~\forall x \in P.\end{align*}  The proof of the theorem is completed now.
 \hfill $\Box$ 
 
 \begin{rmrk}\label{conjugacy} The examples of CCR flows $\alpha^{S^A}$ and $\alpha^{S^{A_T}}$ considered in this Section are in standard form, with vacuum vector providing the invariant states. Further, since they admit only one unit up to scalars, the gauge group action is transitive. So these examples, among themselves, are cocycle conjugate  if and only if  they are conjugate. 
  \end{rmrk}

\section{The Gauge group} 

In this section we describe the gauge group invariant of CCR flows. The computation of the gauge group is used to distinguish CCR flows up to cocycle conjugacy.
   
Recall that a gauge cocycle for  an $\en$-semigroup  $\alpha$ on $B(H)$ is a unitary cocycle $U = \{U_x\}_{x \in P}$ satisfying  $U_x \in \alpha_x(B(H))^\prime$ for all $x\in P$.  
Under the multiplication $(UV)_x:=U_xV_x$, the collection of all gauge cocycles forms a group,  called  the gauge group of $\alpha$. We denote the gauge group by  $G(\alpha)$, 
We endow $G(\alpha)$ with the topology of uniform convergence on compact subsets of $P$ where the the unitary group of $H$ is given the strong operator topology.   
Gauge group is a cocycle conjugacy invariant, whose proof is  
exactly similar to the $1-$parameter case in  \cite[Proposition 2.8.2]{Arveson}. We need a slightly stronger statement of the following Remark, to distinguish our examples.

\begin{rmrk} For any   $\en$-semigroup  $\alpha$ over $P\subseteq \bbr^d$ on $B(H)$, the additive group $\bbr^d$ forms a normal subgroup, identified with 
scalar multiplication by $\{e^{i\langle \lambda | x\rangle }\}_{x \in P}$ for $\lambda \in \bbr^d$.  We denote the quotient group of this normal group in $G(\alpha)$ by $G_0(\alpha)$.  Since these scalars are preserved under the identification between the gauge of two cocycle conjugate $\en$-semigroups, $G_0(\alpha)$ is also a cocycle conjugacy invariant. 
\end{rmrk}

For any  $V:P\to B(K)$, a strictly continuous semigroup of isometric representation, we denote $\m_V = \{V_x, V_x^* :  x \in P \}^{\prime}$. Let $E_x = 1-V_xV_x^*$.  For a unitary $u \in \mathcal{U}(\m_V)$, we define $u_x = u E_x +(1-E_{x})$.

Define $G_V$, as a set, by $G_V = \bbr^d \times \mathscr{A}^V(P) \times \mathcal{U}(\m_V).$   Notice for $h, g \in \mathscr{A}^V(P)$ the continuous function $\lambda(x) = \langle h_x| g_x \rangle$ satisfies $\lambda(x+y) = \lambda(x) +\lambda(y)$ for all $x,y \in P$. So, by Lemma \ref{linext}, there exists a $c(h,g) \in \bbc^d$ satisfying $\langle h_x| g_x \rangle = \langle x | c(h,g)\rangle$ for all $x \in P$.   
Now define multiplication on $G_V$ by $$(\lambda, h, u) (\mu, g, v) = (\lambda+ \mu - Im(c(h, ug)),  h+ug, uv),$$ for $\lambda, \mu \in \bbr^d, h, g \in \mathscr{A}^V(P), u, v \in \mathcal{U}(\m_V)$. It is easy to verify that $G_V$ forms a group under this multiplication. We endow $G_V$ with the product topology given by the standard  topology on $\bbr^d$, the topology of uniform convergence on compact subsets of $P$ on $\mathscr{A}^V(P)$ and the strong operator topology on $\mathcal{U}(\m)$.

\begin{thm} \label{CCR gauge}
 Let $\alpha^V$  be the CCR flow  associated to a pure  strictly continuous semigroup of isometric representation $V$. Then $G(\alpha^V)$ is isomorphic to $G_V$, with the isomorphism given by 
 \[(\lambda, h, u) \longrightarrow \{e^{i\langle \lambda | x\rangle } W(h_x) \Gamma(u_x)\}_{ x \in P}.\] 
 \end{thm}

 \textit{Proof.}  Given any $(\lambda, h, u) \in G_V$ it is easy to verify that $\{e^{i\langle \lambda | x\rangle } W(h_x) \Gamma(u_x): x \in P\}$ forms a gauge cocycle for $\alpha^V$, and that the product structure coincides under the above identification.   If $(\lambda^{(n)}, h^{(n)}, u^{(n)}) \rightarrow (\lambda, h, u)$ in the topology of $G_V$, then it is not difficult to verify that  $e^{i\langle \lambda^{(n)} | x\rangle } W(h^{(n)}_x) \Gamma(u^{(n)}_x)$ strong converges uniformly on compact subsets of $P$, by using the continuity of the exponential vectors.  Conversely if  $U^{(n)}_x=e^{i\langle \lambda^{(n)} | x\rangle } W(h^{(n)}_x) \Gamma(u^{(n)}_x)$ converges uniformly on compact subsets of $P$ in the strong operator topology, then we may conclude as before, by projecting $U^{(n)}_x e(0)$, onto $\bbc e(0)$ and onto the $1-$particle space, that $h^{(n)}$ converges uniformly on compact subsets of $K$ to $h$.  Now the convergence of $U^{(n)}_x e(z)$ uniformly on compact subsets of $P$ imply, again by projecting onto the $1-$particle space, that $u^{(n)} E_x z$  converges to $ u E_x z$ uniformly over compact subsets of $P$. Now the set $\{E_{x}z: x \in P, z \in K\} $  is total in $K$, 
 it follows that $u^{(n)}$ converges to $u$ strongly. We have verified that the topologies of $G_V$ and $G(\alpha^V)$ coincide.	
 
 It can also be verified, by evaluating on exponential vectors, that the map is injective.  
To prove surjectivity, let $\{U_x: x  \in P\} \in G(\alpha)$. For any $a\in \Omega$, the one parameter semigroup of isometries $V^a:=\{V_{ta}\}_{t \in \bbr_+}$  is pure. Thanks to \cite[Theorem 3.8.4]{Arveson}, there exist $\lambda_a, \in \bbr$ , $\{ h^a_{t}\}_{t \in \bbr_+}$ an additive cocycle  for $V^a$  
 and a unitary $u^a \in \{V_{ta} : t \in \bbr_+ \}'$ such that $U_{ta}=e^{it\lambda_a}W(h^a_t)\Gamma(u^aE_{ta}+1-E_{ta}).$ 
 Fix $a,b\in \Omega$. Using the fact that $\{U_{x}\}_{x \in P}$ is a unitary cocycle, compute as follows to observe that 
 \begin{align}
     U_{ta+sb}&=U_{ta}\alpha_{ta}^{V}\left(U_{sb}\right)\cr
     &=e^{it\lambda_a} W(h^a_t) \Gamma(u^{a}E_{ta}+1-E_{ta} ) e^{is\lambda_b}W(V_{ta} h^b_s) \Gamma(E_{ta} +  V_{ta} (u^bE_{sb}+1-E_{sb})V_{ta}^* )\cr
     &=e^{i(t\lambda_a+s\lambda_b)}  W(h^a_t+ V_{ta}h^b_s)  \Gamma\left(u^aE_{ta} + V_{ta}(u^bE_{sb}+1-E_{sb})V_{ta}^*\right)  
                                                             \label{e:1}
       \end{align}
       In the above computation we have used $E_{ta}V_{ta}=0$. Similarly on the other hand, we have
\begin{align}
 U_{ta+sb}&=U_{sb}\alpha_{sb}^{V}\left(U_{ta}\right)\cr    
     &= e^{i(t\lambda_a+s\lambda_b)}  W(h^b_s+ V_{sb}h^a_t)  \Gamma\left(u^bE_{sb} + V_{sb}( u^aE_{ta}+(1-E_{ta})) V_{sb}^*\right) 
                                                            \label{e:2}
     \end{align}

 From Equations \ref{e:1} and \ref{e:2}, we have  
\begin{eqnarray}\label{multiplicative}
u^aE_{ta} + V_{ta}(u^bE_{sb}+(1-E_{sb})V_{ta}^* = u^bE_{sb} + V_{sb}(u^aE_{ta}+1-E_{ta}) V_{sb}^* 
\end{eqnarray} 
and \begin{eqnarray} \label{additive}
     h^a_t+ V_{ta}h^b_s=h^b_s+ V_{sb}h^a_t.
    \end{eqnarray}
For $a\in \Omega$. Set $h_a=h_1^a$, from \ref{additive} we can write   \[ h_a+ V_{a}h_b=h_b+ V_{b}h_a\;\; \text{for every}\; a,b\in \Omega.\]
By letting
$s\to \infty$ in \ref{multiplicative}, we conclude that
 \be \label{uniqueU}
 u^aE_{ta} + V_{ta} u^b V_{ta}^* = u^b  \ee Now letting $t\to \infty $ we get $u^a =u^b$.
Also, by applying  $V_{ta}$ in right on both sides of Equation \ref{uniqueU} (and setting $t=1$), we have $V_{a}u^b=u^b V_{a}$. Since $u^b$ is unitary, it commutes with $V_a^*$ as well. 
Now, since $a \in  \Omega$ is arbitrary and $\Omega $ is dense in $P$,  by the strong continuity, we have $u_b \in \m_V$. 
Thus we have shown that there exists a $u \in \mathcal{U}(\m_V)$ such that for every $a \in \Omega$, there exists $\phi(a) \in \mathbb{T}$, $h_{a} \in K$ such that $V_{a}^{*}h_{a}=0$ and
\[
U_{a}=\phi(a)W(h_a)\Gamma(u_{a}).\]
The fact that the map $\mathbb{T} \times K \times \mathcal{U}(K) \ni (\lambda,\xi,U) \to \lambda W(\xi)\Gamma(U) \in \mathcal{U}(\Gamma(K))$ is injective implies that  $\Omega \ni \to a \to \phi(a) \in \mathbb{T}$ and $\Omega \ni a \to h_{a} \in K$ are well defined maps. That $\phi$ is a homomorphism and $\{h_{a}\}_{a \in \Omega}$ is an additive cocycle follows from the cocycle condition of $\{U_{a}\}_{a \in \Omega}$. Now Lemma \ref{linext} implies that there exists $\lambda \in \mathbb{R}^{d}$ such that $\phi(a)=e^{i \langle \lambda|a \rangle}$. By Lemma \ref{Paddcocycle}, $h$ extends to an additive cocycle on $P$, which we still denote by $h$. Now the gauge cocycle $\{e^{i\langle x, \lambda\rangle}W(h_x)\Gamma(u_x ): x\in P\}$ coincides with $U$ on the dense subset $\Omega$ hence is same as $U$.  This completes the proof. \hfill $\Box$

\begin{rmrk}
\label{commutant of shift}
If $u \in B(L^2(\bbr_+^d)\otimes \k)$ commutes with $\{S_x: x \in P\}$, then $u$ is of the form  $I_{L^2(\bbr_+^d)} \otimes u$ for some $u \in \mathcal{U}(\k)$. This  implies that the commutant of the von Neumann algebra generated by $\{S_{x}:s \in P\}$ is $\{1 \otimes T: T \in B(\k)\}$. In the proof of the following corollary and in the next Section we use this  well-known fact. 
\end{rmrk}

The following corollary, which follows from Theorem \ref{standard form}, Theorem \ref{CCR gauge} and above Remark,  provides a countably infinite family of  mutually non-cocycle-conjugate  
$\en$-semigroups over $\bbr_+^d$.

\begin{cor}\label{c1} Let $\alpha^\k$ be the CCR flow associated with usual  shift of multiplicity $\dim{\k}$, as discussed in Example \ref{usualshift}. Then $G_0(\alpha^\k)$  is isomorphic to $\mathcal{U}(\k)$. 
 
 In particular $\alpha^{\k_1}$ is cocycle conjugate to $\alpha^{\k_2}$ if and only if $\alpha^{\k_1}$ is conjugate to $\alpha^{\k_2}$  if and only if $\dim(\k_1) = \dim(\k_1)$.
\end{cor}

\section{Uncountable many 2-parameter CCR flows}
We end our article by showing that there exists uncountable many examples of CCR flows which are pairwise non-cocycle conjugate, when $P=\mathbb{R}_{+}^{2}$. We need a bit of preprations before we can achieve this.

Let $\ H _1, \ H _2$ be Hilbert spaces, and let $V^{(1)}, V^{(2)}$ be  isometric representations of $P$ on $\ H _1$ and $\ H _2$ respectively. Suppose $T:\ H _2 \to \ H _1$ is a bounded linear Operator, we say that $T$ intertwines $V^{(2)}$ and $V^{(1)}$ if for every $x\in P$   \[V_x^{(1)}T=TV_x^{(2)}\;\; \text{ and} \;\; V_x^{(1)*}T=TV_x^{(2)*}.\]
We denote the space of all bounded linear operators which intertwines $V^{(2)}$ and $V^{(1)}$ by $L(V^{(2)}, V^{(1)})$.  Notice, when $V^{(1)} =V^{(2)}=V$,  the intertwinner  space  $L(V, V)$  is same as the $\m_V$ defined in the previous Section.

Let $a\in \mathbb{R}$ be given. Consider the Hilbert space $L^2(a,\infty)$. For $x\in \mathbb{R}_+$, define $S_x^{(a)}:L^2(a,\infty) \to L^2(a,\infty)$ by 
\begin{equation*}
S_{x}^{(a)}(f)(y):=\begin{cases}
 f(y-x)  & \mbox{ if
} y -x \in [a,\infty) ,\cr
    0 &  \mbox{ if } y-x \notin [a,\infty) .
         \end{cases}
\end{equation*}
Then the map $[0,\infty) \ni x \mapsto S_x^{(a)}\in B\left( L^2(a,\infty)\right)$ is an isometric representation of $\mathbb{R}_{+}$ acting on $L^2(a,\infty).$ 
Fix $a,b\in \mathbb{R}$.  The $1$-parameter isometric representations $\{S_{x}^{(a)}\}_{x\in \mathbb{R}_+}$, $\{S_{x}^{(b)}\}_{x\in \mathbb{R}_+}$ are unitarily equivalent. 
If $U_{b,a}:L^2(a,\infty)\to L^2(b,\infty)$ be the unitary defined by 
$
U_{b,a}f(y)=f(y-b+a)$, then $U_{b,a} S_{x}^{(a)}  =S_{x}^{(b)} U_{b,a}$ for all $x \in \mathbb{R}_{+}$.

For Borel subsets $A,B$ of $ \mathbb{R}$, denote the projection of $L^2(B)$ onto $L^2(B\cap A)$ by $P^B_A$. We simply write $P^B_A$ by $P_A$ when $B$ is clear from the context. 
Let $a_1,a_2\in (-\infty,0)$ be given. Define \[A_1:=[1, \infty)\times [a_1,0) \cup [0,\infty )\times [0,\infty),\]
\[A_2:=[1, \infty)\times [a_2,0) \cup [0,\infty )\times [0,\infty).\]
Note that $A_1,\;A_2$ are $\mathbb{R}_+^2$-modules. 
With respect to the notation introduced in Example \ref{A12} we observe that
$ A_1=A_T^{(-1,-a_1)}+(1,a_1)$ and $A_2=A_T^{(-1,-a_2)}+(1,a_2).$
For $i\in \{1,2\}$, we denote the isometric representation of $\mathbb{R}_+^2$ associated to $A_i$ of multiplicity $1$  
by $V^{(i)}$.

\begin{ppsn} 
If $a_1\neq a_2$, then
  $L(V^{(2)}, V^{(1)})=\{0\}.$ 
\end{ppsn}

$\textit{Proof.} $  Without loss of generality, we can assume that $a_2<a_1$. 
Let  $T\in L(V^{(2)}, V^{(1)})$ be given.
Then for every $(x,y)\in \mathbb{R}^2_+,$ we have 
\[  V^{(1)}_{(x,y)}T=TV^{(2)}_{(x,y)}\;\; \text{and} \;\;  V^{(1)*}_{(x,y)}T=TV^{(2)*}_{(x,y)} . \]
\\
Now write  $A_2$ and $A_1$ as 
 \[A_2:=[1, \infty)\times [a_2,\infty) \cup [0,1 )\times [0,\infty)\;\; \text{ and}  \;\;A_1:=[1, \infty)\times [a_1,\infty) \cup [0,1 )\times [0,\infty).\] 
 Decompose $L^2(A_2)$ and $L^2(A_1)$ as
 \begin{equation} \label{p1ds1}
  L^2(A_2)= K_1\oplus  K_2\;\; \text{and }\;\; L^2(A_1)=\widetilde K_1\oplus \widetilde K_2,
 \end{equation}
 where 
 $ K_1 = L^2(1,\infty) \otimes L^2(a_2,\infty), \;\;  \widetilde K_1 =L^2(1,\infty) \otimes L^2(a_1,\infty), \;\; 
 \text{and} \;\;
       K_2= \widetilde K_2 =L^2(0,1) \otimes L^2(0,\infty)  \quad \text{respectively}.$
Write $T$ in the block matrix form with respect to the decomposition \ref{p1ds1} as
$
   T=(T_{ij})_{2 \times 2}$.
Observe that with respect to this decomposition \ref{p1ds1}, we have
\begin{equation*} 
V_{(0,y)}^{(2)} =\begin{pmatrix} 
        1 \otimes S^{(a_2)}_y  & 0 
  \\0 & 1 \otimes S^{(0)}_y 
  \end{pmatrix} \;\; \text{and }\; \; V_{(0,y)}^{(1)} =\begin{pmatrix} 
        1 \otimes S^{(a_1)}_y  & 0 
  \\0 & 1 \otimes S^{(0)}_y 
  \end{pmatrix} .
  \end{equation*}
The fact  $ V^{(1)}_{(0,y)}T=TV^{(2)}_{(0,y)}$ and $ V^{(1)*}_{(0,y)}T=TV^{(2)*}_{(0,y)}$  for every $y\in \mathbb{R}_+ $ implies that
\begin{equation} \label{p1co1} 
 \begin{pmatrix} 
        (1 \otimes S^{(a_1)}_y)T_{11}  & (1 \otimes S^{(a_1)}_y)T_{12}  
  \\ (1 \otimes S^{(0)}_y)T_{21} & (1 \otimes S^{(0)}_y) T_{22}
 \end{pmatrix}
  = \begin{pmatrix} 
       T_{11} (1 \otimes S^{(a_2)}_y)  & T_{12}(1 \otimes S^{(0)}_y)
  \\ T_{21} (1 \otimes S^{(a_2)}_y) & T_{22} (1 \otimes S^{(0)}_y) 
\end{pmatrix},
  \end{equation}
and
  \begin{equation} \label{p1co2} 
 \begin{pmatrix} 
        (1 \otimes S^{(a_1)*}_y)T_{11}  & (1 \otimes S^{(a_1)*}_y)T_{12}  
  \\ (1 \otimes S^{(0)*}_y)T_{21} & (1 \otimes S^{(0)*}_y) T_{22}
 \end{pmatrix}
  = \begin{pmatrix} 
       T_{11} (1 \otimes S^{(a_2)*}_y)  & T_{12}(1 \otimes S^{(0)*}_y)
  \\ T_{21} (1 \otimes S^{(a_2)*}_y) & T_{22} (1 \otimes S^{(0)*}_y) 
\end{pmatrix}.
  \end{equation}
By equating the $(1,1)^{th}$ entries of Eq. \ref{p1co1} and Eq. \ref{p1co2}, we conclude that $T_{11}(1\otimes U_{a_2,a_1})$ lies in the commutant of the von Neumann algebra generated by $\{1\otimes S_y^{(a_1)} : \; y\in \mathbb{R}_+\}$. By Remark \ref{commutant of shift}, $T_{11}$ is of the form $X_{11}\otimes U_{a_2,a_1}^*$. By comparing other entries of \ref{p1co1} and \ref{p1co2} we conclude that $T$ is of the form   
  \[ T=(T_{ij})=  \begin{pmatrix} 
        X_{11} \otimes U_{a_2,a_1}^*  & X_{12} \otimes U_{a_1,0}  
  \\ X_{21} \otimes U_{a_2,0}^* & X_{22} \otimes 1
 \end{pmatrix}.\]
 Decompose $L^2(a_1,\infty)$ and $L^2(a_2,\infty)$ as $L^2(a_1,0)\oplus L^2(0,\infty)$ and $L^2(a_2,a_1)\oplus L^2(a_1,0)\oplus L^2(0,\infty)$ respectively. Then with respect to the decomposition 
 \begin{center}
 $ L^2(A_2)=\left( L^2(1,\infty)\otimes L^2(a_2,a_1)\right) \oplus \left( L^2(1,\infty)\otimes L^2(a_1,0)\right) \oplus \left( L^2(1,\infty)\otimes L^2(0,\infty) \right) \oplus \left( L^2(0,1)\otimes L^2(0,\infty)\right),$
 \end{center}
 \begin{center}
    $ L^2(A_1)=\left( L^2(1,\infty)\otimes L^2(a_1,0)\right) \oplus \left( L^2(1,\infty)\otimes L^2(0,\infty)\right) \oplus\left( L^2(0,1)\otimes L^2(0,\infty)\right),$
 \end{center}
  T is of the form
 \footnotesize 
 \begin{equation}   \label{p1finaleqn1} 
 T = \begin{pmatrix}
          X_{11} \otimes P_{(a_1,0)}U_{a_2,a_1}^{*}P_{(a_2,a_1)}  &  X_{11} \otimes P_{(a_1,0)}U_{a_2,a_1}^{*}P_{(a_1,0)}      &   X_{11} \otimes P_{(a_1,0)}U_{a_2,a_1}^{*}P_{(0,\infty)}     &   X_{12} \otimes P_{(a_1,0)}U_{a_1,0} 
  \\   X_{11} \otimes P_{(0,\infty)}U_{a_2,a_1}^{*}P_{(a_2,a_1)}  &   X_{11} \otimes P_{(0,\infty)}U_{a_2,a_1}^{*}P_{(a_1,0)}      &    X_{11} \otimes P_{(0,\infty)}U_{a_2,a_1}^{*}P_{(0,\infty)}     &   X_{12} \otimes P_{(0,\infty)}U_{a_1,0} 
  \\  X_{21} \otimes U_{a_2,0}^{*}P_{(a_2,a_1)}   &  X_{21} \otimes U_{a_2,0}^{*}P_{(a_1,0)}     &   X_{21} \otimes U_{a_2,0}^{*}P_{(0,\infty)}    &  X_{22} \otimes 1 
\end{pmatrix}
\end{equation}
\normalsize 
Proceeding in a similar way with the equalities 
$ V^{(1)}_{(x,0)}T=TV^{(2)}_{(x,0)}$ and $V^{(1)*}_{(x,0)}T=TV^{(2)*}_{(x,0)}$ for all $x\in \mathbb{R}_+ ,$
we obtain that $T$ is of the form
  \begin{equation}
  T=
      \begin{pmatrix}  \label{p1finaleqn2}
   1 \otimes  Y_{11}  & 1 \otimes  Y_{12} & U_{1,0}P_{(1,\infty)} \otimes  Y_{13} & U_{1,0}P_{(0,1)} \otimes  Y_{13}
  \\P_{(1,\infty)} U_{1,0}^{*} \otimes  Y_{21} & P_{(1,\infty)}U_{1,0}^{*} \otimes  Y_{22} & 1 \otimes  Y_{23} & 0
\\P_{(0,1)} U_{1,0}^{*} \otimes  Y_{21} & P_{(0,1)}U_{1,0}^{*} \otimes  Y_{22}& 0 & 1 \otimes  Y_{23} 
  \end{pmatrix}
 \end{equation}
 Since $(2,4)^{th}$ and $(3,3)^{th}$ entries of \ref{p1finaleqn2} are zero, the corresponding entries in \ref{p1finaleqn1} are zero which implies that $X_{12}=X_{21}=0.$ Hence $(1,4)^{th}$, $(3,1)^{th}$ and $(3,2)^{th}$ entries of \ref{p1finaleqn1} are zero. Comparing their entries with \ref{p1finaleqn2}, we obtain $Y_{13}=Y_{21}=Y_{22}=0$. This implies that the $(2,2)^{th}$-entry of \ref{p1finaleqn2} is zero. But $P_{(0,\infty)}U_{a_2,a_1}^{*}P_{(a_1,0)} \neq 0$. Thus $X_{11}=0$.  This in turn implies that $Y_{23}=0$. As a consequence, we have $T=0$. This completes the proof. 

\hfill $\Box$

Now we assume that $a_1=a_2 =a$, $A_1=A_2 =A$ and $V^{(1)}=V^{(2)}=V$.

\begin{ppsn} \label{p2}
$L(V, V)= \m_V=  \mathbb{C}1.$ 
\end{ppsn}

$\textit{Proof.} $  
Let  $T\in \m_V$,
then for every $(x,y)\in \mathbb{R}^2_+,$ we have 
$ V_{(x,y)}T=TV_{(x,y)}$ and  $V_{(x,y)}^*T=TV_{(x,y)}^*$.
As before write $A$ as 
$A:=[1, \infty)\times [a,\infty) \cup [0,1 )\times [0,\infty),$ and  
 decompose $L^2(A)$ as
 \begin{equation} \label{ds1}
  L^2(A)= K_1\oplus  K_2,
 \end{equation}
 where 
 $    K_1 =L^2(1,\infty) \otimes L^2(a,\infty), \;\; 
 \text{and} \;\;
        K_2 =L^2(0,1) \otimes L^2(0,\infty)  \quad \text{respectively}.$

Write $T$ and $V_{(0,y)}$ in the block matrix form with respect to the decomposition \ref{ds1} as
\begin{equation*}
   T=(T_{ij})_{2 \times 2} \;\; \text{and}  \;\; V_{(0,y)} =\begin{pmatrix} 
        1 \otimes S^{(a)}_y  & 0 
  \\0 & 1 \otimes S^{(0)}_y 
  \end{pmatrix} . 
  \end{equation*}
The fact $ V_{(0,y)}T=TV_{(0,y)}$ and $ V_{(0,y)}^*T=TV_{(0,y)}^*$  for every $y\in \mathbb{R}_+ $ implies that
\begin{equation} \label{co1} 
 \begin{pmatrix} 
        (1 \otimes S^{(a)}_y)T_{11}  & (1 \otimes S^{(a)}_y)T_{12}  
  \\ (1 \otimes S^{(0)}_y)T_{21} & (1 \otimes S^{(0)}_y) T_{22}
 \end{pmatrix}
  = \begin{pmatrix} 
       T_{11} (1 \otimes S^{(a)}_y)  & T_{12}(1 \otimes S^{(0)}_y)
  \\ T_{21} (1 \otimes S^{(a)}_y) & T_{22} (1 \otimes S^{(0)}_y) 
\end{pmatrix},
  \end{equation}
and
  \begin{equation} \label{co2} 
 \begin{pmatrix} 
        (1 \otimes S^{(a)*}_y)T_{11}  & (1 \otimes S^{(a)*}_y)T_{12}  
  \\ (1 \otimes S^{(0)*}_y)T_{21} & (1 \otimes S^{(0)*}_y) T_{22}
 \end{pmatrix}
  = \begin{pmatrix} 
       T_{11} (1 \otimes S^{(a)*}_y)  & T_{12}(1 \otimes S^{(0)*}_y)
  \\ T_{21} (1 \otimes S^{(a)*}_y) & T_{22} (1 \otimes S^{(0)*}_y) 
\end{pmatrix}.
  \end{equation}
By equating the $(1,2)^{th}$ entries of \ref{co1} and \ref{co2} we conclude that $(1\otimes U_{a,0})^*T_{12}$ lies in the commutant of the von Neumann algebra generated by $\{1\otimes S_y^{(0)} : \; y\in \mathbb{R}_+\}$. By Remark \ref{commutant of shift}, $T_{11}$ is of the form $X_{11}\otimes U_{a,0}$. By comparing other entries of \ref{co1} and \ref{co2} we conclude that $T$ is of the form   
  \[ T=(T_{ij})=  \begin{pmatrix} 
        X_{11} \otimes 1  & X_{12} \otimes U_{a,0}  
  \\ X_{21} \otimes U_{a,0}^* & X_{22} \otimes 1
 \end{pmatrix}.\]
 Now decompose $L^2(a,\infty)$ as $L^2(a,0)\oplus L^2(0,\infty)$. Then with respect to the decomposition 
 \begin{eqnarray*}
 L^2(A)&=& L^2(1,\infty)\otimes L^2(a,0)\oplus L^2(1,\infty)\otimes L^2(0,\infty) \oplus L^2(0,1)\otimes L^2(0,\infty),
 \end{eqnarray*} 
  T is of the form
 \begin{equation}   \label{finaleqn1} 
 T = \begin{pmatrix}
          X_{11} \otimes 1  &    0     &   X_{12} \otimes P_{(a,0)}U_{a,0} 
  \\   0  &   X_{11} \otimes 1        &   X_{12} \otimes P_{(0,\infty)}U_{a,0} 
  \\  X_{21} \otimes U_{a,0}^{*}P_{(a,0)}     &   X_{21} \otimes U_{a,0}^{*}P_{(0,\infty)}    &  X_{22} \otimes 1 
\end{pmatrix}
\end{equation}
 
Proceeding in a similar way with the equalities 
$V_{(x,0)}T=TV_{(x,0)}$ and $V_{(x,0)}^*T=TV_{(x,0)}^*$ for all $x\in \mathbb{R}_+,$
we obtain that $T$ is of the form
  \begin{equation}
  T=
      \begin{pmatrix}  \label{finaleqn2}
   1 \otimes  Y_{11}  & U_{1,0}P_{(1,\infty)} \otimes  Y_{12} & U_{1,0}P_{(0,1)} \otimes  Y_{12}
  \\P_{(1,\infty)} U_{1,0}^{*} \otimes  Y_{21} & 1 \otimes  Y_{22} & 0
\\P_{(0,1)} U_{1,0}^{*} \otimes  Y_{21} & 0 & 1 \otimes  Y_{22} 
  \end{pmatrix}
 \end{equation}
comparing \ref{finaleqn1} and \ref{finaleqn2}, we conclude that $T=\lambda \;1$, for some $\lambda \in \mathbb{C}.$
   \hfill $\Box$

Let $\k_1, \k_2$ be separable Hilbert spaces with $\dim(\k_1), \dim(\k_2)  \in \mathbb{N}\cup \{\infty\}$. Also let  $V^{(A_i,\k_i)}$ be the isometric representation of $\mathbb{R}_+^2$ associated to $A_i$ of multiplicity $\dim(\k_i)$, and $\alpha^{(A_i,\k_i)}$ be the CCR-flow associated to isometric representation $V^{(A_i,\k_i)}$, for $i\in \{1,2\}$. From the above discussions and from Proposition \ref{0addits} the following Corollary is immediate.

\begin{cor} \label{r}
 For $i =1,2$, we have the following.
\begin{enumerate}
 \item [(i)]    If $a_1\neq a_2$ then $L\left(V^{(A_2,\k_2)},V^{(A_1,\k_1)}\right)=\{0\}.$   
\item[(ii)]    The  von Neumann algebra $\m_{V^{(A_i,\k_i)}}$ is isomorphic to $B(\k_i)$.
\item[(iii)]    $G_0(\alpha^{V^{(A_i,\k_i)}})$ is isomorphic to  $\mathcal{U}(\k_i)$.
\end{enumerate}
\end{cor}

The following theorem exhibits uncountably many non-cocycle-conjugate CCR flows.

\begin{thm}  Let $\alpha^{(A_1,\k_1)} $ and $\alpha^{(A_2,\k_2)} $ be the CCR-flows associated to isometric representations $V^{(A_1,\k_1)}$ and $V^{(A_2,\k_2)}$ respectively. Then following statements are equivalent 

\begin{enumerate}
 \item [(i)]    The isometric representations  $V^{(A_1,\k_1)}$ and $V^{(A_2,\k_2)}$  are conjugate. 
 \item[(iii)]    The $\en-$semigroups $\alpha^{(A1,\k_1)} $ and  $\alpha^{(A2,\k_2)}$ are conjugate
\item[(iii)]    The $\en-$semigroups $\alpha^{(A1,\k_1)} $ and  $\alpha^{(A2,\k_2)}$ are cocycle conjugate
\item[(iv)]    $\dim(\k_1)=\dim(\k_2)$ and  $a_1=a_2$.
\end{enumerate}
\end{thm}

\textit{Proof.}  The only non-trivial implication we need to prove  is $(iii) \implies (iv)$.  So we assume that $\alpha^{(A1,\k_1)}$ is cocycle conjugate to
$\alpha^{(A2,\k_2)}$. Since the Lie group dimension  $\mathcal{U}(\k)$ is $\dim(\k)^2$, we immediately conclude, from (iii) of Corollary \ref{r}$,  that\dim(\k_1)=\dim(\k_2)$. So we assume $\k_1=\k_2=\k$.
 


Now suppose that $a_1 \neq a_2$. 
Our assumption that  $\alpha^{(A_1,\k_1)} $ is cocycle conjugate to $\alpha^{(A_2,\k_2)} $, implies that $\alpha^{(A_1,\k)} \otimes \alpha^{(A_1,\bbc)}$ ($\cong \alpha^{(A_1,\k\oplus \bbc)}$) is cocycle conjugate to $\alpha^{(A_2,\k)}\otimes \alpha^{(A_1,\bbc)} .$ 
Note that 
$\alpha^{(A_2,\k)}\otimes \alpha^{(A_1,\bbc)}$ is the CCR-flow associated to the isometric representation  $V^{(A_2,\k)}\oplus  V^{(A_1,\bbc)}.$ By Corollary \ref{r}, we obtain
\begin{eqnarray*} \m_{V^{(A_2,\k)}\oplus V^{(A_1,\bbc)}} & = L\left(V^{(A_2,\k)}\oplus V^{(A_1,\bbc)},V^{(A_2,\k)}\oplus V^{(A_1,\bbc)} \right)\\ & =L\left(V^{(A_2,\k)},V^{(A_2,\k)}\right) \oplus L\left(V^{(A_1,\bbc)},V^{(A_1,\bbc)}\right).\end{eqnarray*}
Hence, by Corollary \ref{r} and Theorem \ref{CCR gauge}, $G_0(\alpha^{(A_2,\k)}\otimes \alpha^{(A_1,\bbc)})$ is isomorphic to $\mathcal{U}(\k)\times \mathbb{T}.$
But $G_0(\alpha^{(A_1,\k\oplus \bbc)})$ is $\mathcal{U}(\k\oplus \bbc)$, which is clearly not isomorphic to  $\mathcal{U}(\k)\times \mathbb{T}$. 
This is a contradiction. Hence $a_1=a_2$. This completes the proof. \hfill $\Box$


\end{document}